\input amstex
\documentstyle{amsppt}
%\magnification=\magstep1
%\pagewidth {12.5 cm}
%\pageheight {19 cm}
\topmatter
\title A $d$-shifted Darboux theorem
%Definitions and conjectures.
\endtitle
\author  E. Bouaziz, I. Grojnowski \endauthor
\address DPMMS, Centre for Mathematical Sciences, 
Wilberforce Road, Cambridge CB3 0WB
\endaddress
%\subjclass Primary 20G99
\abstract{
We give a %etale 
local model for $d$-shifted symplectic dg-schemes, or Deligne-Mumford 
dg-stacks \cite{PTVV}. 
%Etale locally any such is a ``twisted shifted cotangent bundle'', 
Locally any such is a product of a ``twisted shifted cotangent bundle'', 
where the twist is 
given by an element $df$, with $f \in H^{1-d}(\Cal O)$,
and a quadratic bundle in middle degree. 
The latter only occurs if $d = 4r+2$. %This follows from the following.    

\comment
The cotangent complex of a $d$-shifted symplectic stack defines a 
class in the Witt group of the stack; by definition this class is zero
when the cotangent complex is metabolic, that is when there is a 
Lagrangian subbundle. We prove that when this happens the subbundle can be
exponentiated--- the stack is a twisted shifted cotangent bundle.
\endcomment

%The proof is a variant on the classical proofs; but somewhat easier as the
%$d$-shift means there are almost canonical Lagrangian foliations.
} \endabstract

%\endsubjclass
\NoBlackBoxes
\date February 22, 2013. Preliminary version\footnote{Many typos. Some imprecisions. Bad prose.}\enddate
\NoBlackBoxes
\endtopmatter
\document

\define\md#1#2#3#4{\left(\smallmatrix #1 & #2 \\
 #3 & #4 \endsmallmatrix\right)}

\define\rem#1{\medskip\noindent{\bf #1}}

\define\bz{\bold Z}

\define\bA{{\bar{ A}}}
\define\bAh{{\hat{\bar{ A}}}}
\define\bE{{\bar{ E}}}
\define\bD{{\bar{ D}}}
\define\bom{{\bar\omega}}

\define\co{\bold c}

\define\CL{\Cal L}
\define\CV{\Cal V}
\define\CO{\Cal O}
\define\CP{\Cal P}
\define\CA{{\Cal A}}
\define\CX{{\Cal X}}
\define\CY{{\Cal Y}}

\define\LL{{\bold L}}
\define\Lie{{\CL}}

\define\txi{{\tilde\xi}}
\define\teta{{\tilde\eta}}

\define\Mod{\operatorname{mod}}
\redefine\mod{\text{-mod}}
\define\ch{Ch}
\define\Hom{\operatorname{Hom}}

\define\Tor{\operatorname{Tor}}

\define\colim{\operatorname{colim}}

\define\Sym{\operatorname{Sym}}
\define\Spec{\operatorname{Spec}}
\define\tordim{\operatorname{TorDim}}
\define\Disc{\text Disc}
\define\Filt{\operatorname{Filt}}
\define\bFilt{\operatorname{{''}Filt}}

\define\oddline{{k^{0 \mid 1}}}

\define\cone{\operatorname{cone}}
\define\Perf{\text {Perf}\,}

\define\Kos{{\text {Kos}}}

\define\cdga{{\bold {cdga}}}
\define\cdgal{{\bold { cdga_{\leq 0}}}}

\define\Coh{Coh\,}
\define\Pic{Pic\,}

\define\gr{gr}
\define\bgr{\,''gr}
\define\std{{std}}

\define\grade#1{{\langle #1 \rangle}}
\define\intpart#1{\lfloor #1 \rfloor}
\define\mintpart#1{\lceil #1 \rceil}

\redefine\sec#1{\medpagebreak\noindent {#1} \quad}
\define\lsec#1{{\noindent \rm {#1}} \enskip}

\define\quis{{@>{\sim}>>}}

\define\en#1{{\tilde{#1}}}

Let $\CX$ be a derived algebraic variety, or more generally a derived stack \cite{TV,L1,Pr}, defined over a field $k$ of characteristic zero,
and let $\LL_\CX$ be its cotangent complex. Let $\CP \in \Pic \Coh \CO_\CX$ be
an invertible complex of coherent $\Cal D_\CX$-modules, equipped with a 
quasi-isomorphism $\lambda_\CP : (\CP^*)^{-1} \to \CP$.

A {\it $\CP$-shifted symplectic structure on $\CX$} is a $\CP$-valued
deRham closed $2$-form
%$$ \omega = \omega_2 + h \omega_3 + h^2 \omega_4 + \dots \in  \CP \otimes 
$$ \omega = \omega_2 + \omega_3 + \omega_4 + \dots \in  %H^2(
\CP \otimes F^2\LL\Omega_\CX $$
such that 
%i) $(d+D)\omega = 0$,
%
i) $\omega_2 : \LL_{\CX}^* \to \LL_\CX \otimes \CP $ is a quis, and
ii) the induced map $w_2^* \otimes \lambda_\CP : \LL_\CX^* \to \LL_\CX \otimes \CP^*{}^{-1} \to \LL_\CX \otimes \CP $
coincides with $\omega_2$.

%\noindent
When $\CP = \CO_\CX[-d]$,
and for appropriate choice of sign $\lambda_{\CO_\CX[-d]}$,
this is precisely the $d$-shifted symplectic structure of 
\cite{PTVV}: see also \cite{Co}, \cite{AKSZ} and the extensive physics 
literature for earlier avatars of this
definition.

\smallskip

A $\CP$-shifted symplectic %symplectic 
structure on $\CX$ defines an element $[\omega_2]$  
in  the Grothendieck-Witt group $KH_0(\CX,\CP,\lambda_\CP)$, the $0$'th
Hermitian $K$-groups of $\CX$ \cite{Sc}, refining the class of $\LL_\CX $ in 
$K_0(\CX)$.

\smallskip

%In this paper we prove a local structure theorem for the
% $d$-shifted structures, generalizing the classical Darboux theorem, which is
%the  case $d=0$, and the case $d=1$ \cite{BBDJ}.

There is a nice class of examples of $\CP$-shifted symplectic structures,
the {\it twisted shifted cotangent bundles}. We outline their 
definition, repeating it in \S3.2 in more detail.
%\footnote{The details are in the body of the paper in the language of cdga's.}

Let $\CY$ be a derived algebraic variety such that $\LL_\CY$ is
perfect, let $\CP = \CO_\CY[-d]$, for $d > 0$, $\lambda_\CP = \pm 1$,
and suppose 
the sheaf $H^i(\CP^{-1} \otimes \LL_\CY^* ) = 0$ for $i > 0$.

Then the shifted cotangent bundle $\bar\CX = \LL_\CY \otimes \CP
= \Spec \Sym_{\CO_\CY}(\CP^{-1} \otimes \LL_{\CY}^*  ) \to \CY$
canonically has a $\CP$-shifted symplectic structure $d\bar\lambda$,
where $\bar\lambda : \CO_{\bar\CX} \to \CP \otimes \LL_\CX $ is the Liouville
$1$-form, defined just as in classical symplectic geometry.

\comment
Here we have made a slight departure from usual derived algebraic geometry,
and consider everything enriched not in the category of chain complexes
of vector spaces, but in the category of chain complexes of 
 $\bz/2\bz$-graded vector spaces.
This convention allows us to happily ignore the signs $\lambda_\CP$
(the symplectic monoidal structure takes care of them), %  for us),
and treat symmetric %($\lambda_\CP =1 $) 
and anti-symmetric %($\lambda_\CP = -1$)
forms on equal footing, just by replacing $\CP$ by  $\CP \otimes \oddline$ 
if $\lambda_\CP = 1$, where  $\oddline$ denotes the one dimensional 
vector space $k$ with parity $-1$, the odd line.  
\endcomment

 Now let
$\txi \in H^0(\CY,H^2(\CP \otimes F^1\LL\Omega_\CY))$, so 
%$\txi = \txi_1 +h\txi_2 + \cdots $ is de Rham closed, with 
$\txi = \txi_1 +\txi_2 + \cdots $ is de Rham closed, with 
$\txi_1 : \CO_\CY[1] \to \CP \otimes \LL_\CY$.

Then $\txi_1$ defines a one parameter family of deformations of $\bar\CX$,
whose generic fiber is $\CX =  
%(\bold A^1_\CY[1] @>{\xi_1}>>  \LL_\CY \otimes (\CP')^* ) \times_{ Y \times \bold A^1} Y =
 \Spec \Sym_{\CO_\CY}(\CP^{-1} \otimes \LL_{\CY}^* @>{\xi_1}>> \CO_\CY[-1]) 
\times_{ Y \times \bold A^1} Y @>{\pi}>> Y$,
and the deformation $\lambda$ of the Liouville form $\bar\lambda$ to $\CX$ 
combines with $\pi^*\txi$ to define $\omega_\txi : = (d+D)\lambda - \pi^* \txi
\in H^2(\CX,\CP \otimes F^2\LL\Omega_\CX)$, a $\CP$-shifted symplectic structure.
\smallskip
Given a $\CP$-shifted symplectic stack $\CX$, and $\CV \in \Perf \CX$ a 
complex with $H^i(\CV) = 0 $ for $i \geq 0$ and 
$\alpha : \CV^\dagger \to \CV$ a quis such that 
$\alpha  = \lambda_\CP \alpha^\dagger$,  
where $\CV^\dagger = \CV^* \otimes {\CP}^{-1}$,
 then 
$\CX \times \CV^* = \Spec \Sym_{\CO_\CX} \CV$ is also naturally 
a $\CP$-shifted symplectic stack.

 If $d \neq 2 \Mod 4$ this will not produce
any new examples locally, but for $d = 2 \Mod 4$ gives us new examples
attached to non-zero classes in the Witt group of quadratic forms.
%This obstruction can  always 
%be trivialised etale locally, and can mostly be trivialised Zariski locally.
%However if  $\CP = \CO[-d]$, with $d = 4r + 2> 0$ an integer congruent to
%$2$ mod $4$, a

%There is another class of examples of $\CP$-shifted sympectic structures

\smallskip

In this paper we prove that locally every $\CP$-shifted symplectic structure
on $\CX$ is  a product of  a {twisted shifted cotangent bundle}
and a quadratic bundle $\CV$, with $\CV$ zero unless $d =2\Mod 4$.
This
generalize both the classical Darboux theorem, which is
the  case $d=0$, and the theorem of \cite{BBDJ},\footnote{
As we were %in the final stages of 
writing this note, Brav, Brussi and Joyce released a preprint, arXiv/1305.6302, %with a similar result.
with a different normal form for $-d$-symplectic varieties. They give a Hamiltonian description. A twisted cotangent bundle has such a description: if  
$\CY = \Spec \Sym k[z_1,\dots,z_n | Dz_i = h_i ]$, $\xi = df$, $f \in H^{1-d}(\CO_\CY)$, and $y_i = dz_i^\dagger$, then $H = f + \sum Dz_i.y_i$.}
 which is part of
the case $d=1$.

\smallskip

We remark that for $d = 2 \Mod 4$, any quadratic bundle
on a variety $Y$ with a non-trivial class in the Witt group of the function
field $k(Y)$ of $Y$ defines a $\CP$-shifted symplectic variety $\CX$ 
on the shifted 
cotangent bundle of $Y$. 
If the dimension of the quadratic bundle is even,
it is etale locally, but 
{\it not} Zariski locally, a metabolic bundle, and so $\CX$ is still  etale locally
a twisted shifted cotangent bundle,
and if the dimension of the bundle is odd we can write it etale locally
as the product of a one dimension quadratic bundle with a metabolic bundle.

\smallskip

The proof has 4 steps.
We begin by assuming $\CX$ is foliated by almost Lagrangian subvarieties---that
is, there is a smooth morphism 
$ \pi: \CX \to \CY$
such that the vertical maps in
$$\CD 
\pi^* \LL_\CY  @>>> &  \LL_\CX   @>>> &  \LL_{\CX/\CY} @>{+1}>> \\
  @A{\omega_2}AA    & @A{\omega_2}AA   & @A{\omega_2}AA \\
\LL_{\CX/\CY}^\dagger @>>> & \LL_\CX^\dagger  @>>> & \pi^* \LL_\CY^\dagger @>{+1}>>
\endCD$$
are quasi-isomorphisms, where $\CA ^\dagger = \CA^* \otimes {\CP}^{-1}$.

This is not necessarily possible---there is an 
obstruction to doing so if the class 
of $\omega_2$ in the Witt group of $\CX$ is non-zero. 

In step 1 of the proof we show that the class of $\omega_2$ in the Witt
group of $\CX$ is the only obstruction to %locally 
finding such a foliation.
The main ingredient for this is a dg-Frobenius theorem, 
stating when a `subbundle' of the cotangent complex can be integrated to
a map of dg-schemes, 
which we formulate and
prove in proposition 1.4.

In step 2, we show that, locally on $\CX$,
the map $\CX \to \CY$ degenerates in a dg-smooth
family to $\bar{\CX} = T^*[d]\CY \to \CY$, the $d$-shifted cotangent 
bundle of $\CY$.

Unlike smooth families of ordinary schemes, this does not imply the de Rham
cohomology is constant in the family (any semi-stable %log smooth 
family is dg-smooth--- for
example a smooth elliptic curve degenerating to a nodal one). However
we can always define a specialisation map for closed $p$-forms
which become %specialize to 
zero in the de Rham complex
of the generic fibre; we get a closed form on the special fibre.\footnote{This is only an issue when $d = 1$. When $d>1$, we can choose $\CY$
so $\pi_0(\CX) \to \pi_0(\CY)$ is an isomorphism, 
and we have 
$H^i(\CX,\LL\Omega_\CX) \quis H^i(\CY,\LL\Omega_\CY ) \quis H^i(\bar{\CX},\LL\Omega_{\bar{\CX}}) $. In this case the composite can be described explicitly on the level of chain complexes in terms of flat sections of the Gauss-Manin connection of the degenerating family.}

Using this, we can transport the class $\omega \in F^2\LL\Omega_\CX$;
it becomes a class %$\psi(\omega) \in
in $ F^1\LL\Omega_{\bar{\CX}}$,
which is the pullback of a class in $F^1\LL\Omega_{\bar{\CY}}$.
This class defines a deformation of $\bar{\CX}$ to $\CX'$, a twisted
shifted cotangent bundle.

In step 3 we 
% use the Moser technique\cite{Mo} to
show
that there is an automorphism of $\CP \otimes \LL_\CY$ which induces
an isomorphism from $\CX$ to $\CX'$ sending $\omega $ to the standard
form on $\CX'$. 
This can be accomplished by the Moser technique; % \cite{Mo}; 
%we're embarrassed to admit 
our setup is already sufficiently simple that we just do it by hand. 
%
% given two $\CP$ shifted symmetric structures $\omega,\omega^\std$
%on $\CX$, if $H^{2}(\CX,\CP\otimes\LL\Omega_\CX) = 0$,
%there is an automorphism of $\CX$ which takes $\omega$
% to $\omega^\std$.
%
%This proves the theorem.

Finally, in step 4 we observe that if the class of $\omega_2$ is non-zero
in the Witt group, then over a (Zariski!) local neighbourhood we can write 
$\CX$ as a product of a quadratic bundle and a 
$\CP$-shifted symplectic variety with $\omega_2$ zero in the Witt group.
This is an immediate consequence of Witt cancellation and algebraic surgery.
% from the theory of derived quadratic forms.

\smallskip 
%\footnote{
The proofs are written in the body of the paper, perhaps perversely, 
without using geometric language. 

%Finally, in step 3 we show %use the Moser technique\cite{Mo} to show that given a $\CP$ shifted symmetric structure $\omega$ on $T^*\CY^\txi$, if $H^{2-d}(\CY,\LL\Omega_\CY) = 0$, there is an automorphism of $\CY$ which takes $\omega$ to the standard symplectic structure, proving the theorem.

\medskip

%It is clear that the above sketch of proof, and its precise implementation
%below extend without much effort to mixed and non-zero characteristic
%when we replace deRham with crystalline cohomology. This will be presented elsewhere 

\medskip
The contents of \S0 and \S1 are taken from notes for a course one of us taught in 2005, and make no particular claim to originality. \S2 is more background, 
much of which is unneccesary to the proof, but which we find pedagogically
reassuring. Finally, in \S3 we prove the theorem.

%Both groups have been aware 

\comment
This note is part of a larger project, which includes \cite{FG}.
\endcomment

%We first encountered examples of $d$-symplectic varieties in our work \cite{FG}.

%We first encountered examples of $d$-symplectic varieties in our work \cite{FG} on the Hilbert scheme.

%This paper is part of a project to understand some examples of 

%

\head Conventions \endhead
In a model category, we write $\quis, \hookrightarrow, \twoheadrightarrow$ to mean morphisms which are, respectively,  weak equivalences, cofibrations or fibrations. %, and denote $\co B \to B$ a cofibrant replacement of $B$.

We sometimes write $A/B$, and sometimes  $\cone (B \to A)$, for the mapping cone
of a morphism $B \to A$.

Morphisms are strict unless otherwise indicated.

\head 0. The cotangent complex \endhead

The various (Quillen equivalent) modern formulations of
$(\infty,1)$-categories give a clean conceptual picture of the
cotangent complex and what it means. This formidable technology is
unnecessary %irrelevent
 to the theorems of this paper, which are about local
computations, and for which the straightforward and classical theory
of model categories suffice. We recall briefly the %necessary
definitions \cite{Qu,BG,GS}.

\sec{0.1} We fix a field $k$, $char(k) =0$. Denote by 
$Ch(k)$ the category of chain complexes of $k$-modules, 
%$Ch(k)$ the category of chain complexes of $\bz/2\bz$-graded $k$-modules, 
$\dots\to M^i @>{\delta^i}>> M^{i+1} \to \dots$,
equipped with the projective model category structure and the usual symmetric
monoidal
structure.
Weak equivalences are quasi-isomorphisms, and fibrations are degree-wise surjections. Let $\cdga$ denote the category of commutative differential graded algebras over $k$, the commutative algebra objects in $\ch(k)$. The adjunction
$$ Sym : Ch(k) \leftrightarrows \cdga : Forget $$
allows us to transport the model category structure on $\ch(k)$ to one on $\cdga$
%\cite{BG},
and the adjunction 
$$ Include : \cdgal \leftrightarrows \cdga : \tau_{\leq 0} $$
allows us to transport the model category structure on $\cdga$ to $\cdgal$, the
category of non-positively graded cdga's over $k$.
Explicitly, a morphism $f : A \to B $ in $\cdgal$ is a weak equivalence if it
is a quis, that is $H^n(f) : H^n(A) \to H^n(B)$ is an isomorphism for all $n$, 
 $f$ is a fibration if $f^n : A^n \to B^n$ is a surjection for all $n \leq -1$,
and $f$ is a cofibration if it is one in $\cdga$. These have been described in 
\cite{BG}, and we recall:
% that 
a set of generating cofibrations are  given by  $k \to \Sym k[n]$,
$k \to \Sym \Disc_n$, and $\Sym k[n] \to \Sym \Disc_n$, for $n \leq 0$, where 
$\Disc_n = k[n-1] @>{Id}>> k[n]$.
Cofibrations are closed under  pushouts, a retract of a cofibration is a cofibration, if $X_1\to X_2 \to \dots $ is a sequence of cofibrations, then $X_1\to \colim X_i $ is a cofibration, and if $f_j : X_j \to Y_j$, $j \in J$ is a set of cofibrations, then $\otimes f_j : \otimes X_j \to \otimes Y_j$ is a cofibration. 

%If $A \in \cdgal$ and $\alpha \in A^n$ is an $n$-cycle, we write $A[x \mid Dx = \alpha]$ for the pushout $A \otimes_{\Sym k[n]} \Disc_n$.

% and so a typical cofibration is of the form $B \to B[z_\alpha 

\bigskip
The model category $\cdga$, and hence $\cdgal$, is left proper; that is if $f : A \to B$ is a quis,
and $A \hookrightarrow C$ is a cofibration, then $C \to B \otimes^L_A C$ is a quis.

\bigskip
If $A \in \cdga$, $A\mod$ is the dg-category whose objects are the objects $M \in \ch(k)$ endowed with a strictly associative action of $A$. The free $A$-module functor gives an adjunction
$$ A \otimes_k (\cdot) :  Ch(k) \leftrightarrows A\mod : Forget, $$
the induced model category structure on $A\mod$ is called the projective model category, and is enriched over $\ch(k)$. More generally, given a cofibration %{\it check?}
$A \to B$ in $\cdga$, the adjunction $ B \otimes_A (\cdot) :  A\mod \leftrightarrows B\mod : Forget $ is Quillen when both $A\mod$ and $B\mod$ are equipped with the projective model category structure, and is a Quillen equivalence if $A\to B$ is a quis.

\bigskip

An object $B$ in a model category $\Cal C$ enriched in chain complexes
is said to be {\it compact} if $\Hom_{\Cal C}(B,\cdot)$ commutes with 
filtered homotopy colimits. If $\Cal C$ has a fixed set of compact generators, 
we say $B$ is {\it finitely presented}, or f.p.{}, if it is in the smallest subcategory containing the generators and closed under shifts, finite coproducts, mapping cones and weak equivalences. For such a category an object is compact if and only if it is a retract of a f.p.{} object. We apply these notions to $\cdgal$, $\cdga_{A \backslash \cdot}$, and to $A\mod$, where we denote the category of compact objects $\Perf A$.
% and note they are precisely the dualisable complexes, 
%Observe if $A \to B$ is a quis, then Perf A and Perf B are Quillen equivalent

\sec{0.2} If $f : A \to B$ is in $\cdga$, we write $\LL_{B/A}$ for the
cotangent complex. Recall the  definition: Take a cofibrant replacement 
for $f$ in the model category $\cdga_{A\backslash\cdot}$ of alegbras over $A$,
i.e. factor $f = pi$,   $A @>{i}>> \co B @>{p}>> B$, with $i $ a cofibration and $p$ a weak equivalence. Then $\LL_{B/A} = B \otimes_{\co B} \Omega^1_{\co B/A} \in B\mod$.

As $\cdga$ is left proper, this coincides with the cofibrant replacement in the
model category of morphisms in $\cdga$. In other words, if we first 
cofibrantly replace $A$, $\co A \quis A$, and then factor $\co A \to B$ by
$ \co A \hookrightarrow X \quis B$, the natural map 
$B \otimes_X \Omega^1_{X/\co A} \to \LL_{B/A}$ is a quis.

\bigskip
We recall  (i) a square $\CD A & @>>> &  R \\  @VVV & & @VVV \\ B & @>>> & S \endCD$
induces a morphism
 $S \otimes_B \LL_{B/A} \to \LL_{S/R}$ in $S\mod$, 
that this is a weak equivalence if both $ B \to S$, $A \to R$ are, and

ii) morphisms $A \to B \to C$ in $\cdga$ induce a triangle
$$ C \otimes_B \LL_{B/A} \to \LL_{C/A} \to \LL_{C/B} @>{+1}>> $$
in $C\mod$.

\sec{0.3} Given $f : A \to B$ in $\cdga$, factor $ f = pi : A @>{i}>> \co B
 @>{p}>> B$ as in 0.2, and define
$$ F^p \LL\Omega_{B/A} = \prod_{n \geq p} {\bigwedge}^n_{\co B} \Omega^1_{\co B/A}[-n]
\qquad \in A\mod,$$
with differential $D+d$, where $D$ is induced from the differential $D$
on $\Omega_{\co B/A}^1$, and $d$ is the derivation induced from 
$d : \co B \to \Omega^1_{\co B/A}$.

%{\it Fuck I'm tired. We could also take the 

Write $[ \,\, ]_p : F^p\LL\Omega_{B/A} \to \bigwedge^p\LL_{B/A}  = 
B \otimes_{\co B} \bigwedge^p_{\co B}\Omega^1_{\co B/A} $ for the morphism
%of $B$-modules 
induced by the projection and a shift by $p$.

%These complexes are well defined, up to weak equivalence.

\smallskip
We write $\LL\Omega_{B/A}$ instead of $F^0\LL\Omega_{B/A}$, and call
this the Hodge completed de Rham complex \cite{FT}.
If $A = k$, we write $\LL_B$, $F^p\LL\Omega_B$ instead 
of $\LL_{B/k}$, $F^p\LL\Omega_{B/k}$. We also write 
$H^i(\Spec B, F^p\LL\Omega_B)$ for the $i$'th cohomology of the 
complex $F^p\LL\Omega_{B/k}$.

% Hodge completed de Rham complex

\sec{0.4} Recall \cite{Go} that if $f : B \to A $ in $\cdgal$ is such that
$H^0(B) \to H^0(A)$ has nilpotent kernel, then the induced map
$\LL\Omega_B \to \LL\Omega_A$ is a quis in $k\mod$.

In particular, if $B \in \cdgal$, then $H^i(\Spec B, \LL\Omega_B) = 0$
for $i <0$.

\head 1.Basics % what to call this
\endhead

\sec{1.1}
We say that a morphism $X \to Y$ is $k$-connected if $H^i\cone(X \to Y) = 0$ for all $ i \geq -k$.\footnote{This is probably usually called $-k$-co-connected. Whatever.} Observe $X \to Y$ is $k$-connected if and only if 
$ H^{-k} X \to H^{-k}Y$ is an epimorphism and,  for $i > -k$, 
$ H^i (X) \to H^i(Y) $ is an isomorphism.

\proclaim{Proposition} Let $B \in \cdgal$ be f.p.{}, and $B \to A$ a morphism in $\cdga$. Then for $d \geq 1$ the following are equivalent.

i) $H^0B \to H^0A$ is a surjection, %isomorphism, 
and $H^i(\LL_{A/B}) = 0$ for $i \geq -d +1$,
and

ii) $H^iB \to H^iA$ is an isomorphism for all $i > -d + 1$,
and $H^{-d+1}B \to H^{-d+1}A$ is a surjection, i.e.{} $B \to A$ is $(d-1)$-connected.

\smallskip

\noindent
Moreover, if $K = \cone(B \to A) \in B\mod$, then 
 $H^{-d}(K) = H^{-d}(A \otimes^L_B K) 
\quis H^{-d}( \LL_{ A/B})$.

\endproclaim
\demo{Proof} Replace $ B \to A$ by a cofibration, $ B \hookrightarrow \co A \quis A$. 
The small object argument shows that $\co A$ is built out of $B$ by 
a (possibly transfinite) attachment of cells; that is there is an ordinal $I$, and for each $\alpha \in I$ an element $f_\alpha \in B[z_\beta , \, \beta < \alpha]$, such that 
$ \co A := B[z_\alpha,\, \alpha \in I \mid Dz_\alpha = f_\alpha]$.\footnote{
If $A \in \cdgal$ and $\alpha \in A^n$ is an $n$-cycle, we write 
$A[x \mid Dx = \alpha]$ for the pushout $A \otimes_{\Sym k[n]} \Disc_n$,
and we say $\deg x = n -1$.
}
Furthermore, if  $H^iB \to H^iA$ is an isomorphism for all $i > -r + 1$,
and $H^{-r+1}B \to H^{-r+1}A$ is a surjection,
then we can insist $\deg z_\alpha \leq -r$.

We have $\LL_{A/B} = A\langle dz_\alpha, \, \alpha \in I \mid D(dz_\alpha) = -df_\alpha \rangle$.

Let $K $ be the mapping cone of $B \to \co A$ in $B\mod$. Then $K = B\{ z_\alpha \}$, where, given elements $\gamma_\alpha \in \co A$, we write $B\{ \gamma_\alpha\} $ for the $B$-subalgebra of $\co A$ they generate. %d by the elements $\gamma_\alpha$.
The derivation $ d :\co A  \to \LL_{\co A/B} $ factors through $K$,  and the mapping cone of the induced map $1\otimes d : \co A \otimes_B K \to \LL_{\co A/B}$ in $\co A\mod$ is $Q := \co A\{1\otimes z_\alpha z_\beta 
%- z_\alpha \otimes z_\beta  - (-1)^{|z_\alpha||z_\beta|} z_\beta \otimes z_\alpha 
- z_\alpha \otimes z_\beta  - (-1)^{\deg z_\alpha \deg z_\beta} z_\beta \otimes z_\alpha 
, \, \alpha,\beta \in I \} [1] $.
Hence if $\deg z_\alpha \leq -r$ for all $ \alpha  \in I$, then $ H^iQ = 0 $ for $ i > - 2r -1 $.

Now, $H^iK = 0 $ for $i > - r$ if and only if $H^i(\co A \otimes_B K) = 0 $ for $i > - r$, as is evident from the spectral sequence 
$Tor^{HB}(HA, HK) % =  HA \otimes^L_{HB} HK
 \Rightarrow H(\co A \otimes_B K)$.

Fix $ r $ such that $ H^iK = 0 $ when $ i \geq - r + 1$ and $ r \geq 1$; 
this is possible as $H^0B \to H^0A$ is a surjection. % isomorphism. 
Then $H^{i-1}Q = H^iQ = 0$ for $i \geq - 2r$, $H^i(\co A \otimes_B K) = 0 $ for $i \geq -r + 1$, and the triangle $ \co A \otimes_B K \to \LL_{\co A/B} \to Q @>{+1}>> $ gives $H^i (\LL_{\co A/B}) = 0 $
for $i \geq -r +1 $ and an isomorphism $H^{-r}(K) = H^{-r}(\co A \otimes_B K) 
\quis H^{-r}( \LL_{\co A/B})$.

Hence (i) implies (ii), on taking $r = d$, but also (ii) implies (i).

\enddemo

\proclaim{Corollary} Let $f : B \to A$ in $\cdgal$  have $\LL_{A/B} = 0$.

(i) If $H^0B \to H^0A$ is an
isomorphism, then $f $ is a quis. %quasi-isomorphism.

(ii)  If $H^0B \to H^0A$ is a
surjection, then $f $ is a quis. %quasi-isomorphism.

\endproclaim
\comment
\demo{Proof} 
(i) is immediate. For (ii), observe that with notation as in
the proposition, $H^0(B) \to H^0(A)$ a surjection implies we can insist
$\deg z_\alpha \leq -1$, and that $H^i(K) = 0$ when $i \geq 0$. Taking
$ r= 1$ gives an isomorphism $H^{-1}(K) \quis H^{-1}(\LL_{\co A/B}) = 0$,
and so $H^0B \to H^0A$ is an isomorphism. The result follows from (i).
\enddemo
\endcomment

\sec{1.2}
We say that $M \in A\mod$ has Tor amplitude in $[a,b]$
if for all
$N \in A\mod$ such that $H^j(N)=0$ for $j \neq 0$, $ H^i(M \otimes_A N) \neq 0$
implies $a \leq i \leq b$, and that it has  Tor dimension $d$, $\tordim M = d$, 
if it has Tor amplitude in $[-d,0]$.

%Equivalently, [if M Perfect] it suffices to insist that 
% for all $k$-points $\varphi : A\to k$, $H^i(M \otimes_A k) \neq 0$ unless $-t \leq i \leq 0$.
% Nakayama's lemma?

Observe $A[d]$ has Tor dimension $d$ if  $d \geq 0$, 
and that if $ M' \to M \to M'' @>{+1}>> $
is a triangle, then $\tordim M \leq \tordim M' + \tordim M''$.
Hence if $M \in \Perf A$, then $M$ has finite Tor dimension.

\proclaim{Lemma} If $M \in A\mod$ has $H^iM = 0$ for $i \geq -r$ and 
if $\tordim M \leq r$, then $M = 0$.
\endproclaim
\demo{Proof} Suppose $M \neq 0$, and let $j$ be the maximum integer such
that $H^j(M) \neq 0$, so $j < -r$. Then the spectral sequence 
$ \Tor^{H(A)}(H(M), H^0(A)) \Rightarrow H(M \otimes^L_A H^0(A)) $
gives $H^j(M \otimes_A H^0(A)) \neq 0$, so $-r \leq j$, a contradiction.
\enddemo
%{\it careful! finiteness needed?)

\sec{1.3}
A morphism $A\to B \in \cdgal$ is a {\it Zariski open embedding} if it is quis in $\cdga_{A\backslash \cdot}$
to a morphism of the form $A \to A[t,\xi | D\xi = tf-1]$ where $f \in A^0$. A {\it cover} of a f.p.{} algebra $A$ is a finite set $f_i : A \to B_i$ of Zariski open embeddings such that if $B = \oplus B_i$, the augmented bar complex $A \to B \to B\otimes_A B \to \cdots$ is exact; equivalently the maps $H^0(f_i) : H^0(A) \to H^0(B_i)$ form a cover in the usual sense. We say a property of an $M \in A\mod$ holds {\it Zariski locally} on $A$ if there is a cover of $A$ for which the the property holds for the pullback $B_i\otimes^L_A M$ on $B_i$ for all $i$.

\proclaim{Lemma} Let $A \in \cdgal$, $H^0(A)$ finitely presented, $M \in \Perf A$. The following are equivalent.

i) $M$ has Tor amplitude in $[a,b]$.

ii) $M \quis (M_a \to \dots \to M_b)$, where each $M_i[-i]$ is a summand of a free
module $A^{n_i}$, for some $n_i \geq 0$.

iii) There is a Zariski cover $A \to A'$ such that %if $M' = M\otimes_A A'$
 $M\otimes_A A' \quis (M'_a \to \dots \to M_b')$, where each $M'_i[-i]$ is a free
module ${A'}^{n_i}$, for some $n_i \geq 0$.

\endproclaim
\comment
\demo{Proof}
\enddemo
\endcomment

\sec{1.4} Let $R \to A\in \cdgal$, and suppose given 
%(i) an exact sequence $ S \hookrightarrow \LL_{A/R} \twoheadrightarrow \LL_{A/R}/S$ in
(i) a triangle 
 $ S \to \LL_{A/R} \to \LL_{A/R}/S @>{+1}>>$ in
$A\mod$, with $\LL_{A/R}$ and $\LL_{A/R}/S$ cofibrant, 
and (ii) a map $d : \LL_{A/R}/S \to \wedge^2 (\LL_{A/R}/S)$ in $R\mod$ with $d^2 = 0$
such that
$\CD   \LL_{A/R}  @>{\beta}>>  \LL_{A/R}/S \\
 @V{d}VV  @V{d}VV \\
 \wedge^2 \LL_{A/R} @>{\wedge^2\beta}>>  \wedge^2 (\LL_{A/R}/S) \endCD $
(strictly) commutes.

Call such data {\it foliation data}; observe that if $B \to A$ is a morphism
in $\cdga_{R\backslash\cdot}$,
we get foliation data from $\LL_{B/R}\otimes_B A \to \LL_{A/R} \to \LL_{A/B}$ after 
cofibrantly replacing $B\to A$.

Define $\LL\Omega_{A/S} = \prod_{i \geq 0} \wedge^i (\LL_{A/R}/S) [-i] \in R\mod$ 
with differential $D + d$, where $D$ is the internal differential,
% a complex in $\ch(k)$,
and observe we have a morphism $\LL\Omega_{A/R} \to \LL\Omega_{A/S}$ in $R\mod$.

\proclaim{Proposition} 
 Let $A \in \cdgal$ have $H^0(A)$ finitely presented as a commutative algebra,
and suppose given either

(i) $S \in \Perf A$, and morphism $ S \to \LL_A$, such that,
writing $\LL_A/S = \cone(S \to \LL_A )$,
there is an $s>0$ with  $H^i(\LL_A/S) =  0$ for $i > -s$, and $\tordim S \leq 2s-1$, or

(ii) foliation data $S \to \LL_A \to\LL_A/S @>{+1}>>$ 
with $S \in \Perf A$, as above.

\noindent
Suppose also that
%$H^0(S) \quis H^0(\LL_A)$ is an isomorphism,
$S \to \LL_A$ is $0$-connected.
\footnote{ A lazy assumption,
but we do not need the general case in this paper.} 
% and $H^iS = 0 $ for $i > 0$.

Then there exists a finitely presented $B \in \cdgal$, morphism $B \to A$,
with $H^0(B) \to H^0(A)$ a surjection, % isomorphism, 
and quis $\LL_B\otimes_B A \quis S$ factoring $\LL_B\otimes_B A \to \LL_A$.
In case (ii) we also have a weak equivalence of triangles
%$\CD \LL_B\otimes_B A @>>> \LL_A @>>> \LL_{B/A} \\
% @VVV @VVV @VVV \\
%S @>>> \LL_A @>>> \LL_A/S \endCD $.
$(\LL_B\otimes_B A \to \LL_A \to \LL_{A/B}) @>>>
(S\to  \LL_A \to \LL_A/S) $
such that 
$\CD \LL_{A/B} @>{d}>> \wedge^2 \LL_{A/B} \\
 @VVV @VVV \\
 \LL_{A}/S @>{d}>> \wedge^2 (\LL_{A}/S) \endCD$ (strictly) commutes.

\endproclaim
\demo{Proof} We inductively build algebras $B_r \to A$,
and morphisms $\LL_{B_r} \otimes_{B_r} A \to  S \to \LL_A$ factoring 
 $\LL_{B_r} \otimes_{B_r}A \to  \LL_A$ 
 such that, writing $S/\LL_{B_r}$ for $\cone( \LL_{B_r}\otimes_{B_r}A \to S)$,
i) $B_r$ is finitely presented, with $\tordim \LL_{B_r} \leq r$,
ii) for $r \geq 1$, $H^0 B_r \to H^0 A$ is a surjection, % isomorphism,
iii) $H^i(S/\LL_{B_r}) = 0$ for $i > -r$, 
iv) $S/\LL_{B_r}$ is perfect, and $\tordim S/\LL_{B_r} \leq \tordim S$,
and (v) when $S$ comes from foliation data, $\CD \LL_{A/B_r} @>{d}>> \wedge^2 \LL_{A/B_r} \\
 @VVV @VVV \\
 \LL_{A}/S @>{d}>> \wedge^2 (\LL_{A}/S) \endCD$ (strictly) commutes.

Then for $r > \tordim S + 1$  lemma 1.2 gives $S/\LL_{B_r} = 0$, proving the proposition.

\smallskip

To begin, as $H^0A$ is finitely presented we can choose $x_1, \dots, x_n$ 
in $A^0$ which generate $H^0(A)$, and $\xi_1,\dots,\xi_r$ in $A^{-1}$ such
that the $D\xi_i$ generate the ideal $D(A^{-1})$ in $A^0$. 
Define $B_0 = k[x_1',\dots,x_n']$,  and
 $\LL_{B_0} \otimes_{B_0} A \to S$ by sending $dx'_i$ to 
some choice of a lift of $dx_i \in \LL_A^0$ to $S^0$; this is possible
as $H^0(S)$ surjects onto $H^0(\LL_A)$.
\comment
Define $B_0 = k[x_1,\dots,x_n]$,  choose lifts $f_i \in B_0$ of $D\xi_i$, and 
define$ B_1 = k[x_1,\dots,x_n,\xi_1,\dots,\xi_r \mid D\xi_i = f_i ]$,
equipped with the evident maps $B_0 \to B_1 \to A$. It is clear that
$B_0$, $B_1$ satisfy the inductive requirements above.
As $S$ is perfect and $H^0(\LL_A)$ is a quotient of 
$A\langle dx_1,\dots,dx_n\rangle$, the $H^0(A)$-module
$\ker ( H^0(S) \to H^0(\LL_A))$ is finitely generated;
let $v_{n+1},\dots,v_{n+m}$ be elements of $S^0$ which generate it.
Define $B_1 = k[x'_0,\dots,x'_{n+m}]$, with $\deg x'_i = 0$,
and define $B_1 \to A$ to be the map which sends $x'_i \mapsto x_i$ if
 $i \leq n$, and $x'_i \mapsto 0$ if $ i > n$.
Define $\LL_{B_1} \otimes_{B_1} A \to S$ by sending $dx'_i$ to 
some choice of a lift of $dx_i \in \LL_A^0$ to $S^0$  for $i \leq n$,
and to $v_i$ if $ i > n$. It is clear that $B_1$ satisfies the
inductive requirements (i)-(iv).
\endcomment

Now suppose we have defined $B_r$. If $H^{-r}(S/\LL_{B_r}) = 0$, set $B_{r+1} = B_r$, otherwise %consider the map 
%the map $S \to S/\LL_{B_r}$ defines a quotient space $\text{Coker}( H^{-r}(S)   \to H^{-r}(S/\LL_{B_r}))$.
as $S$ and $\LL_{B_r}$ are perfect and $H^0(A)$ is finitely presented,
  $H^{-r}(S/\LL_{B_r})$
is a finitely generated $H^0(A)$-module;
let $w_1,\dots,w_t$ be elements which generate it.

As in proposition 1,
 let $K = \cone(B_r \to A) \in B_r\mod$, $d : K \to \LL_{A/B_r}$.
Consider the image of $w_i$ in $H^{-r}(\LL_{A/B_r})$. We claim this is in
the image of $ d : H^{-r}K \to H^{-r}(\LL_{A/B_r})$.

Granting this for the moment, let $w_1',\dots, w_t' \in H^{-r}K$ be 
elements for which $dw_i' $ equals the image of $w_i$ in $H^{-r}(\LL_{A/B_r})$.
These define elements in  $H^{-r+1}(B_r)$ under the transgression in 
the exact sequence
$$ \to H^{-r}A \to H^{-r}(K) \to H^{-r+1}B_r \to H^{-r+1}A \to $$
and we let $\gamma_1',\dots,\gamma_t' \in (B_r)^{-r+1}$ 
be lifts of these elements from $H^{-r+1}(B_r)$.

Write $\gamma_1,\dots,\gamma_t$ in $A^{-r+1}$ for the images of the $\gamma_i'$
under the morphism $B_r \to A$. As these elements are zero in $H^{-r+1}A$,
there exist $z_i \in A^{-r}$ such that $Dz_i = \gamma_i$.

We may further assume that $dz_i = dw_i'$ in $H^{-r}(\LL_{A/B_r})$.
To see this, choose a cofibrant replacement $B_r \to \co A$ for $B_r \to A$. As in
proposition 1.1, we can assume $K \subseteq A$, and so $w_i' \in A^{-r}$
has $Dw_i' = \gamma_i' \in B_r^{-r+1} \subseteq A_r^{-r+1}$. Now choose
$z_i = w_i'$.

%Finally, let $v_1,\dots,v_s$ generate the image  $ H^{-r}(S)   \to H^{-r}(S/\LL_{B_r}))$ as an $H^0(A)$-module.

Set
$ B_{r+1} = B_r[z'_i \mid Dz'_i = \gamma'_i]$, with $\deg z'_i = -r$,
and map $B_{r+1} \to A$ by sending $z'_i \mapsto z_i$.

As $d\gamma_i' \in H^{-r+1}(\LL_{B_r})$ are transgressed from 
$w_i \in H^{-r}(S/\LL_{B_r})$, we can find $\tilde{w}_i \in S^{-r}$
lifting $w_i$ with $D\tilde{w}_i $ equal to the image of $ -d\gamma_i'$
in $S^{-r+1}$. Extend the map $\LL_{B_r} \to S$ to $\LL_{B_{r+1}} \to S$ 
by sending $dz_i \mapsto \tilde{w}_i$.

Then $B_{r+1}$ satisfies the inductive conditions (i)-(iv).

It remains to prove the claim. For case (i) of the proposition,
note that if %$S = \LL_A $,
$H^i(\LL_A/S) = 0$ for $i > -s$, then for $ r \leq s$ and $i > -r \geq -s$,
the exact sequence 
$\to H^i(S/\LL_{B_r}) \to H^i(\LL_{A/B_r}) \to H^i(\LL_A/S) \to$
gives $H^i(\LL_{A/B_r}) = 0$.
Hence proposition 1 gives 
%that if $j$ is the maximum integer with $H^j(K) \neq 0$, then 
%$H^j(K) = H^j(A \otimes^L_{B_r} K) = H^j(\LL_{A/B_r})$, so $j =-r$ 
$H^{-r}(K) = H^{-r}(A \otimes^L_{B_r} K) = H^{-r}(\LL_{A/B_r})$, and the claim.
If $ s  < r \leq 2s$, we reduce to case (ii) by constructing foliation data
Zariski locally.

Replace $B_r \to A$ by $B_r \to \co A$, so $\LL_{A/B_r}$ is cofibrant 
and $(\LL_{A/B_r})^i = 0 $ for $ i > -s$. We have the Tor amplitude 
of $S/\LL_{B_r}$ is contained in $[-2s+1,-r]$.

Suppose 
$$S/\LL_{B_r} \quis A\langle \eta_1,\dots,\eta_i \rangle = A \otimes_k M[1],
\tag **$$
where $D \eta_i \in A \langle \eta_1,\dots,\eta_{i-1}\rangle$,
$M[1] = \oplus k \eta_i $, $ -2s + 1 \leq \deg \eta_i \leq -r$,
and $S/\LL_{B_r} \to \LL_{A/B_r}$ is induced by $\xi : M[1] \to \LL_{A/B_r}$.

Then writing  $C =  \cone(S/\LL_{B_r} \to \LL_{A/B_r}) = \LL_{A/B_r} \oplus (A\otimes M)$  
(with differentials),  
$\wedge^2 C = \wedge^2 \LL_{A/B_r} \oplus (\LL_{A/B_r} \otimes M) \oplus \wedge^2M$
(with differentials),  we can define $d : C \to \wedge^2 C$ by
$(\omega, a\otimes \eta_i) \mapsto (d\omega, da \otimes \eta_i,0)$.
This is a morphism of complexes, as 
$d(a \xi(\eta_i)) = da.\xi(\eta_i) + a.d\xi(\eta_i) = da. \xi(\eta_i)$,
as $d\xi(\eta_i) \in (\wedge^2 \LL_{A/B_r})^{\deg \eta_i + 1} = 0$, as 
$-\deg\eta_i - 1 < 2s$.
Moreover, $d ^2 : C \to \wedge^3 C $ is zero, and $d$ defines foliation
data on $S/\LL_{B_r} \to \LL_{A/B_r} \to C$. %DIAGRAM commutes is clearer

To finish the proof, for both cases of the %!! theorem 
proposition we have,
by construction,  foliation data
$S/\LL_{B_r} \to \LL_{A/B_{r}} \to \LL_{A}/S$, and so we have a morphism
 $\LL\Omega_{A/B_r} \to \LL\Omega_{A/S}$. 
Let $Q[1]$ be its mapping cone. We have 
$Q \quis \prod_{i \geq 1} Q_i[-i]$, where each 
$Q_i  =\cone(\wedge^i \LL_{A/B_r} \to \wedge^i (\LL_A/S)) \in A\mod$ 
has a filtration with subquotients 
$\wedge^a (\LL_{A}/S) \otimes \wedge^b (S/\LL_{B_r})$, $a+b = i$, and $b > 0$.

As $H^k(\LL_A/S) = 0$ for $k \geq 0$, and $H^k(S/\LL_{B_r}) = 0$ 
for $k > -r$,
$H^k(\wedge^a (\LL_A/S) \otimes \wedge^b (S/\LL_{B_r})) = 0$
when $k > -br -a$. Hence $H^{1-i-r}(Q_i) = 0$ when $i > 1$,
and any $w \in H^{-r}(S/\LL_{B_r}) = Q_1$ defines a class in $H^{1-r}Q$,
and hence a class in $Im (H^{1-r}(F^1\LL\Omega_{A/B_r}) \to H^{1-r}(\LL\Omega_{A/B_r}))$.

This is the case even when assumption (**) does not hold, as it
always holds Zariski locally, and so there is a Zariski 
cover $\phi : A \to A'$ such that, for $u \in Im (H^{-r}(S/\LL_{B_r}) \to 
H^{-r}(\LL_{A/B_r}))$, $d\phi(u) =0$. But $d\phi(u) = \phi du$,
and  
as the kernel of $\phi  : H^{-r}(\wedge^2 \LL_{A/B_r}) \to 
H^{-r}(\wedge^2\LL_{A'/B_r}) = H^{-r}(\wedge^2 \LL_{A/B_r}) \otimes_{H^0A} H^0{A'}$
is zero, $du = 0$.

 For $ r\geq 1$, %this image
$Im (H^{1-r}(F^1\LL\Omega_{A/B_r}) \to H^{1-r}(\LL\Omega_{A/B_r}))$
is zero, so the triangle 
$F^1\LL\Omega_{A/B_r} \to \LL\Omega_{A/B_r} \to 
\LL\Omega_{A/B_r}/F^1\LL\Omega_{A/B_r} @>{+1}>>$ 
gives that the image of $w$ is of the form $dw'$, for some $w' \in A^{-r}$,
$D w'=0$, as required.

Finally, we must check the inductive condition (v). Again, choose 
cofibrant replacements  $B_r \to \co A$ for $B_r \to  A$, and 
$\co ( S/\LL_{B_r}) \to  S/\LL_{B_r} $ such that $ (\co S/\LL_{B_r})^i = 0$
for $i > -r$. Then $(Q_i)^{-r} = 0$ for $i > -r$, so $dw = 0$ for 
$w \in (\LL_A/S)^{-r}$, ensuring (v) for $B_{r+1}$.

%\footnote{This is truly awful writing. It would be best if we rewrote it completely. Though currently tortuous, it is a much stronger result then any we need;  the easy special cases already suffice for our Darboux theorem. A slightly better way to organise the material would be to include much more systematic development. Show, given $S$, the category $\Cal C$ of $B \to A$ with $\LL_B\otimes_A A \quis S$ has $B,B' \in \Cal C$ implies there is a $C \in \Cal C$ and quis $C \to B$, $C \to B'$. Show case (i) implies  foliation data  exists Zariski locally, and globally.}

\enddemo

\proclaim{Corollary} 
 An algebra $A \in \cdgal$ is f.p.{} if and only if $H^0(A)$ is finitely presented as a commutative algebra and $\LL_A$ is perfect.
\endproclaim
\demo{Proof} The `only if' implication is clear.
For the converse, apply the proposition with $S\to \LL_A$ equal
to $\LL_A @>Id>> \LL_A$. This gives a f.p.{} algebra $B$ with
$H^0(B) \quis H^0(A)$ an isomorphism  and $\LL_{A/B}= 0$. Corollary 1.1 gives $B\to A$ is a quis.
%corollary  implies that $B_r \to A$ is a quis, proving the theorem.

\enddemo

Note that the proof shows that if $\LL_A$ has Tor dimension $d > 0$ then $A$ is
quis to a cdga built by attaching $r$-cells for $r \geq -d$,  $A \quis k[z_1,\dots,z_n \mid Dz_i = f_i]$ for some elements $f_i \in k[z_1,\dots,z_{i-1}]$ with $\deg z_i \geq -d $ for all $i$.

\comment % OLD VERSION

\proclaim{Proposition} 
 An algebra $A \in \cdgal$ is f.p.{} if and only if $H^0(A)$ is finitely presented as a commutative algebra and $\LL_A$ is perfect.
\endproclaim
\demo{Proof} The `only if' implication is clear.
To prove the converse, we inductively build algebras $B_r \to A$ such that
i) $B_r$ is finitely presented, with $\tordim \LL_{B_r} \leq r$,
ii) for $r \geq 1$, $H^0 B_r \to H^0 A$ is an isomorphism,
iii) $H^i(\LL_{A/B_r}) = 0$ for $i > -r$, and 
iv) $\LL_{A/B_r}$ is perfect, and $\tordim \LL_{A/B_r} \leq \tordim \LL_A$.

Then for $r > \tordim \LL_A + 1$ the lemma gives $\LL_{A/B_r} = 0$, and so the
corollary  implies that $B_r \to A$ is a quis, proving the theorem.

To begin, as $H^0A$ is finitely presented we can choose $x_1, \dots, x_n$ 
in $A^0$ which generate $H^0(A)$, and $\xi_1,\dots,\xi_r$ in $A^{-1}$ such
that the $D\xi_i$ generate the ideal $D(A^{-1})$ in $A^0$. Define
$B_0 = k[x_1,\dots,x_n]$,  choose lifts $f_i \in B_0$ of $D\xi_i$, and 
define
$ B_1 = k[x_1,\dots,x_n,\xi_1,\dots,\xi_r \mid D\xi_i = f_i ]$,
equipped with the evident maps $B_0 \to B_1 \to A$. It is clear that
$B_0$, $B_1$ satisfy the inductive requirements above.

Now suppose we have defined $B_r$. If $H^{-r}(\LL_{A/B_r}) = 0$, set $B_{r+1} = B_r$, otherwise let $K = \cone(B_r \to A) \in B_r\mod$, and let $j$ be the maximal
integer such that $H^j(K) \neq 0$. Proposition 1 gives 
$H^j(K) = H^j(A \otimes^L_{B_r} K) = H^j(\LL_{A/B_r})$, so $j =-r$ and we have
an exact sequence 
$$ \to H^{-r}A \to H^{-r}(\LL_{A/B_r}) \to H^{-r+1}B_r \to H^{-r+1}A \to 0.$$
The map $\LL_A \to \LL_{A/B_r}$ defines a subspace $Im H^{-r}(\LL_A)$ in
$H^{-r}(\LL_{A/B_r})$. %, and hence a subspace in $H^{-r+1}(B_r)$. 
As $\LL_A$ is perfect,
this is a finitely generated $H^0(A)$-module. 
Let $w_1,\dots,w_t$ be elements which generate it,
$\gamma_1',\dots,\gamma_t' \in (B_r)^{-r+1}$ 
lifts of their transgressions to $H^{-r+1}(B_r)$,
and $\gamma_1,\dots,\gamma_t$ in $A^{-r+1}$ the images of the $\gamma_i'$
under the morphism $B_r \to A$. As these elements are zero in $H^{-r+1}A$,
there exist $z_i \in A^{-r}$ such that $Dz_i = \gamma_i$.

Set
$ B_{r+1} = B_r[z'_i \mid Dz'_i = \gamma'_i]$, with $\deg z'_i = -r$,
and map $B_{r+1} \to A$ by sending $z'_i \mapsto z_i$.
Then $B_{r+1}$ satisfies the inductive conditions (i)-(iv).

\enddemo

Note that the proof shows that if $\LL_A$ has Tor dimension $d > 0$ then $A$ is
quis to a cdga built by attaching $r$-cells for $r \geq -d$,  $A \quis k[z_1,\dots,z_n \mid Dz_i = f_i]$ for some elements $f_i \in k[z_1,\dots,z_{i-1}]$ with $\deg z_i \geq -d $ for all $i$.

\endcomment

\proclaim{\lsec{1.5} Corollary} Let $A \in \cdgal$ have $H^0(A) $ f.p.{}. and
$\LL_A^*[-d] \quis \LL_A$ for some $d \geq 0$. Then $\tordim \LL_A \leq d$, and $A$ is f.p.{}
\endproclaim
\demo{Proof} We have $\LL_A^{**} \quis \LL_A$, and so $\LL_A$ is perfect. A perfect complex $M$ has Tor amplitude in $[a,b]$ if and only if $M^*$  has Tor amplitude in $[-b,-a]$, so $\LL_A$ has Tor dimension $d$.
\enddemo

\proclaim{\lsec{1.6} Lemma} Suppose $B \to A$ in $\cdgal$ has $H^0(B) \to H^0(A)$ 
an isomorphism. Let $M,N \in \Perf B$. 
If $\phi : M \otimes_B A \to N\otimes_B A$ is a quis,
then there exists a quis $\bar\phi : M \to N$ in $B\mod$. % inducing $\phi$.
\endproclaim
\demo{Proof} Let $M, N$ have Tor 
amplitude in $[a,b]$, and $H^k(M) = H^k(N) = 0$
for $k > j$, and suppose $j$ is minimal. We induct on $j -a  +  b$,
if this is negative $M=N=0$ and the lemma is true. 
Othewise, as $H^0(B)\to H^0(A)$ is an isomorphism, $H^j(M) \to H^j(N)$ 
is also an isomorphism,  and we can choose $\gamma_1,\dots,\gamma_t \in M^j$ 
which generate $H^j(M)$ as an $H^0(B)$-module. Put $\Gamma = \oplus B\gamma_i$,
$\tilde{M} = \cone(\Gamma @>>> M)$, $\tilde{N} = \cone(\Gamma @>{\phi}>> N)$.
Then $\tilde M, \tilde N$ still have Tor amplitude in $[a,b]$, but
$H^k(\tilde M) = H^k(\tilde N) = 0$ for $k \geq j$.

By induction there exists a quis $\tilde\phi : \tilde M \to \tilde N$,
% inducing $\tilde M \otimes_B A \to \tilde N \otimes_B A$, 
and we set
$\bar\phi $ to be the cone of $\tilde\phi[-1] \oplus Id : \cone(
\tilde M [-1] + \Gamma ) \to \cone( \tilde N [-1] + \Gamma)$

\enddemo

\proclaim{\lsec{1.7} Lemma} Let $A \in \cdgal$ be f.p., and suppose $P,Q \in \Perf A$
satisfy $P \otimes^L_A Q \quis A$. Then Zariski locally on $A$ there is a quis
$ P \quis A[d]$, for some integer $d$ (which may depend on the open set, if $H^0(A)$ is not connected)
\endproclaim

\head 2. More background, didactally %First reductions 
\endhead

%Throughout this section we fix $d \geq 1$. Given $B \in \cdgal$ and $ M\in B\mod$, define 
%$M^\dagger = \Hom_B(M,B)[-d] \in B\mod$. Note that if

\sec{2.1}
Given $M \in B\mod$, $\xi \in H^0Hom_B(M,B[1])$, write $\Sym^\xi_BM$ for the
pushout $\Sym_B \widetilde{M} \otimes_{\Sym_B B} B$, where $\widetilde{M} 
= \cone( M[-1] @>{\xi}>> B)$.

\sec{2.2} 
Fix $B \in \cdgal$ f.p.{}, $M \in \Perf B$ with $H^i(M) = 0$ for $i > 0$, and
$\xi \in H^0Hom_B(M,B[1])$ as above. The algebra $\Sym^\xi_BM$ is
filtered, with $B$ in filtration degree 0 and $M$ in filtration degree 1;
the associated graded algebra is $\Sym_B M$.

Write $A = \Sym^\xi_BM$, $\bA = \Sym_BM$.
The grading on $\bA$ is defined by a $\bold G_m$-action;
differentiating we get a map $Lie {\bold G_m} = k \to \LL_\bA^*$. Denote
the image of $1$ by $\bE$, this is the `Euler vector field',
$ \bE \in H^0(\Spec \bA, \LL_\bA^*) = H^0(Hom(\LL_\bA,\bA))$.
If we choose a cofibrant representative for $M$, so 
$M = B\langle y_1,\dots,y_n \mid Dy_i = \sum_{j<i} \mu_{ij}y_j \rangle$,
$\mu_{ij} \in B$, then $\bE = \sum y_i (dy_i)^*$, $D_\bA \bE = 0$.

\sec{2.3}
The grading on $\bA$ induces one on $\LL_\bA$ and on $\LL\Omega_\bA$;
the $i$'th graded piece is the $i$'th eigenspace of $\Lie_\bE$; $\bA$,
$\wedge^r\LL_\bA$ and $\LL\Omega_\bA$ are the sums of their graded pieces.

\proclaim{Lemma} i) If $\lambda \neq 0$, the $\lambda$'th graded piece of 
$H^i(\Spec \bA, F^p\LL\Omega_\bA)$ is the $\lambda$'th graded piece of 
$ \ker ( d : H^{i-p}(\wedge^p\LL_\bA) \to  H^{i-p}(\wedge^{p+1}\LL_\bA) )$.

ii) The $0$'th graded piece of $H^i(\Spec \bA, F^p\LL\Omega_\bA)$ is
$H^i(\Spec B, F^p \LL \Omega_B)$.
\endproclaim
\demo{Proof}
Choose a cofibrant replacement for $M$ as above.
Define
% $D_\Kos : \LL_\bA \to \bA$ by $D_\Kos(\LL_B) = 0$, $D_\Kos(dy_i) = y_i$. Then $D_\Kos D dy_i = - D_K d (\sum_{j < i} \mu_{ij}y_j) = -D D_\Kos dy_i$, 
$D_\Kos$ to be contraction with $\bE$, so
%so $[D, D_\Kos] = 0 : \LL_\bA \to \bA$, and we can define $D_\Kos : \wedge^i \LL_\bA \to \wedge^{i-1} \LL_\bA$ so that $D_\Kos$ is a derivation of cdga's; 
$[D,D_\Kos] = 0$, and 
$[d, D_\Kos] = [d+ D,D_\Kos] = \Lie_\bE$.

Let $\omega = \omega_p + \omega_{>p} \in H^i(\Spec \bA,  F^P\LL\Omega_\bA)$,
$\Lie_\bE \omega = \lambda \omega$. Then
$ [D_\Kos , d + D] \omega_{>p} = \Lie_\bE \omega_{>p} = \lambda \omega_{> p}$, and
$ d\omega_p = -(d + D) \omega_{>p}$, so $\omega_{>p} = D_\Kos(-d\omega_p/\lambda)
+ (d+D)(D_\Kos\omega_{>p}/\lambda)$,
i.e.{} $\omega = \omega_p + D_\Kos(-d\omega_p/\lambda) 
= dD_\Kos(\omega_p/\lambda)
\in H^i(\Spec \bA, F^P\LL\Omega_\bA)$. This shows the inclusion
$\ker ( d : \wedge^p\LL_\bA[-p] \to  \wedge^{p+1}\LL_\bA[-p-1])  \hookrightarrow F^p\LL_\bA$
is a quis on $\lambda$'th graded pieces, $\lambda \neq 0$.

By definition the $0$'th graded piece of $F^p\LL\Omega_\bA$ is $F^p\LL\Omega_B$,
proving (ii).

\enddemo

\sec{2.4} Regard the class $\xi : M \to B[1]$ as a vector field on 
$\bA$ via the map $M^* \otimes_B \bA \to \LL_\bA^*$.

The vector field $\bE \in H^0(\LL_\bA^*)$ does not deform to a vector field
on $A$, and neither does  the cohomology class of $\xi \in H^1(\LL_\bA^*)$. 
Instead, there is an 
$E \in (\LL_A)^0$ and $\xi \in (\LL_A)^{1}$ such that $D_A E = \xi $,
and the specialisation of $E$, $\xi$ to $\bA$ are $\bE$, $\xi$, respectively.
On choosing a cofibrant representative for $M$ as in 2.2, 
$E = \sum y_i(dy_i)^*$, $\xi = \sum \xi_i (dy_i)^*$.

\smallskip
Choose cofibrant replacements for $B$ and $M$, and an
element $\xi \in \Hom(M[-1],B)^0$ whose class is $\xi$. This defines 
cofibrant representatives 
 for $A$ (resp.{} $\wedge \LL_A$)
having the same 
underlying algebra  as that for $\bA$ (resp.{} as $\wedge \LL_\bA$), 
but $D_A = D_\bA + \Lie_\xi$, where 
$\Lie_\xi = [d,\xi] : \wedge^i \LL_A  \to \wedge^i \LL_A$, $i \geq 0$.

\smallskip

Write $\bD = D_\bA$. The operators $d, \bD, \xi, \Lie_\xi, D_\Kos, \Lie_\bE$
act on the underlying %($\bz/2\bz$-graded) 
chain complex of $\LL\Omega_A$,
%$F^P\LL\Omega_A$, 
and these operators are subject to the relations
$[d,D_\Kos] = \Lie_\bE$, $[D_\Kos,\Lie_\xi] = \xi$, $[\Lie_\bE,\Lie_\xi] = -\Lie_\xi$, 
$[D_\Kos,\xi] = [\bD,D_\Kos] = [\bD, d] = 0$.

\sec{2.5} 
The filtration on $A = \Sym^\xi_BM$ induces  an increasing filtration on
$\wedge^r\LL_A$, denoted $\bFilt^\delta\wedge^r \LL_A$.
There is a related filtration that is also %more 
useful. The triangle
$$ \LL_B \otimes_B A \to \LL_A \to \LL_{B/A} \quis M\otimes_B A @>{+1}>> $$
induces an increasing filtration $\Filt^\delta\wedge^r\LL_A$ of
$\wedge^r\LL_A$;    
$\Filt^\delta\wedge^r\LL_A$ is the $A$-submodule of $\wedge^r\LL_A$
spanned by terms $d\gamma_1\dots d\gamma_r$ with at most $\delta $
of the $\gamma_i \in M$.

Define a decreasing exhaustive filtration on $\LL\Omega_A$ by
$$ \Filt^\delta \LL\Omega_A = \prod_{r \geq 0} \Filt^{r-\delta}\wedge^r\LL_A, $$
and on $F^p \LL\Omega_A$ by
 $\Filt^\delta F^p \LL\Omega_A = \prod_{r \geq p} \Filt^{r-\delta}\wedge^r\LL_A = \Filt^\delta\LL\Omega_A \cap F^p\LL\Omega_A. $

%\remark{Remark}
%The filtration on $\LL\Omega_A$ whose $\delta$'th piece
%is $\prod_{r \geq 0} \bFilt^\delta \wedge^r\LL_A$ has no useful properties---it is 
%not exhaustive, and inclusion of the union of all its pieces into 
%$\LL\Omega_A$ is {\it not} a quis.
%\endremark

Define $\gr^\delta F^p\LL\Omega_A = \Filt^\delta F^p\LL\Omega_A / 
 \Filt^{\delta +1}  F^p\LL\Omega_A $. 

\comment
\proclaim{Lemma}
i) $\gr^\delta F^p\LL\Omega_A = \prod_{k \geq \max(p,\delta)}
\wedge^\delta(\LL_B\otimes_B A) \otimes_A \wedge^{k-\delta}(M \otimes_BA)$,
and the induced differentials $D$,$d$ on $\gr^\delta$ are $B$-linear.

ii) If $p > \delta$, the map
$ (\wedge^\delta \LL_B \otimes \wedge^{p-\delta}M, D) \to
(\gr^\delta F^p\LL\Omega_A, d+D)$
is a quis in $B\mod$.

iii) If $p > \delta$, the map
$\Filt^\delta F^p\LL\Omega_A \hookrightarrow \widetilde{\Filt^\delta F^p\LL\Omega_A} 
:= \{ \omega \in F^p\LL\Omega_A \mid [\omega]_p \in \Filt^{p-\delta}\wedge^p \LL_A \}$
is a quis.

\endproclaim
\demo{Proof} (i) is immediate from the definition.

(ii) Note $d=d_{A/B}$ is identically $0$ on $\wedge^i_BM$, so the map claimed
really is a morphism of complexes. Now observe that the filtration on $A$
induces one on $\gr^\delta F^p\LL\Omega_A$, that the map of (ii) is a filtered
morphism, and the associated graded  is
$ \wedge^\delta \LL_B \otimes \wedge^{p-\delta}M \to \bar{\gr}^\delta
=: \prod_{k \geq \max(p,\delta)}
\wedge^\delta(\LL_B \otimes_B \wedge^{k-\delta}M) \otimes_B\bA$. 
It suffices to show that this is a quis.

Define $D_\Kos : M \to \bA = \Sym_BM$ to be the inclusion,
and extend this to $D_\Kos: (\bar{\gr}^\delta)^i \to  (\bar{\gr}^\delta)^{i-1}$
so $D_\Kos$ is a derivation of cdga's, $[D,D_\Kos] = 0$, and 
$[d,D_\Kos] = \bE$. The result now follows as in lemma 2.3.

(iii) Now, %the cone 
$\cone(\Filt^\delta F^p\LL\Omega_A \to  F^p\LL\Omega_A) $
is a filtered complex, whose associated graded is 
$ \oplus_{\delta < t \leq p} \gr^t  F^p\LL\Omega_A \quis
 \oplus_{\delta < t \leq p}  (\wedge^t \LL_B \otimes \wedge^{t-p}M )\otimes_B A$
by (ii), and this is just the associated graded of the filtered complex
$\cone(\Filt^{p-\delta} \wedge^p\LL_A \to  \wedge^p\LL_A)$.

As  $ \widetilde{\Filt^\delta F^p\LL\Omega_A} = 
\cone( F^p\LL\Omega_A \oplus \Filt^{p-\delta} \wedge^p\LL_A [p]
\to \wedge^p\LL_A [p] )$, we get that the associated graded of 
$\cone(\Filt^\delta \to \widetilde{\Filt^\delta})$ is quis to $0$,
and hence the cone is quis to $0$.

\enddemo

\endcomment

\comment
\remark{Remark} 
%{\it sloppy phrasing}
%On choosing a cofibrant replacement for $M$ as above,
%the vector field $\bE$ deforms to $E \in (\LL_A^*)^0$ with $D_A E = \xi$, $ E = \sum y_i(dy_i^*)$,
%In these terms
On choosing a cofibrant replacement for $M$ as above,
 $\bFilt^\delta\wedge^r\LL_A = \ker( \Lie_E(\Lie_E -1) \dots (\Lie_E - \delta +1 ) : \wedge^r \LL_A \to \wedge^{r} \LL_A$, 
 $\Filt^\delta\wedge^r\LL_A = \ker( E^\delta : \wedge^r \LL_A \to \wedge^{r-\delta} \LL_A)$. 
Note $\ker E \subseteq \ker \xi$, re-confirming that these are $D$-stable

\endremark
\endcomment

\sec{2.6}
If $N \in \Perf B$, and $\xi : N \to B $ a morphism in $B\mod$,
define the completion of $B$ at the ideal generated by $N$, 
$$ \hat{B}_N = \lim_{@<<<} B \otimes^L_{\Sym_B N^{\otimes_B n}} B $$
where $B$ is a $\Sym_B N^{\otimes_B n}$-algebra in two ways:
via the augmentation morphism sending $N^{\otimes_B n}$ to zero,
and via the morphism  ${\xi^{\otimes n}}: Sym_B N^{\otimes_B n} @>>> B$.

\smallskip

Observe that $\LL_{\hat{B}_N} = \LL_B \otimes_B \hat{B}_N$.

\sec{2.7}
We use this notation when $N = M[-1]$, $\xi : M[-1] \to B$ as above, 
and denote
the completion $\hat B_{M[-1]})$ by $\hat B$. 
Also write $\bAh = \hat B \otimes_B \bA$.

\proclaim{Lemma} i) There is a well defined morphism of chain complexes
$e^\xi : F^p\LL\Omega_A \to \LL\Omega_\bAh $,
$$ w = w_p + w_{p+1} + \dots \mapsto
e^\xi\omega = \sum_{a \geq 0} 
\left(\sum_{k \geq \max(0,p-a)} \frac{\xi^k \omega_{k+a}}{k!}\right).$$

ii) This descends to a morphism 
$e^\xi : \Filt^\delta F^p\LL\Omega_A \to \Filt^\delta F^\delta\LL\Omega_\bAh $.

iii) If the image of $\omega \in H^i(\Spec A,\Filt^\delta\LL\Omega_A)$ in 
$H^i(\Spec A, \LL\Omega_A)$ is zero, then $e^\xi\omega$ is the image of 
a class in $H^i(\Spec \bA, \Filt^\delta F^\delta\LL\Omega_\bA) $,
with $[e^\xi\omega]_\delta = \frac{\xi^{p-\delta}}{(p-\delta)!}\omega_p$.

\endproclaim
\demo{Proof} $\xi^k\omega_{k+a}$ is manifestly zero in 
$\bA \otimes_{Sym( M[-1]^{\otimes l})} \bA$ if $ k \geq l$, so the sum
$\sum_{k}  \frac{\xi^k}{k!} \omega_{k+a}$ is a well defined
element of the completion $\wedge^a\LL_\bAh = \wedge^a(\LL_\bA \otimes_\bA \bAh)$.
Moreover $(d + D_A - \Lie_\xi) e^\xi = e^\xi (d + D_A)$,
as $\Lie_\xi = [d,\xi]$ commutes with $\xi$.\footnote{This is a special case of the 
Campbell-Baker-Hausdorff formula, or a straightforward computation in 
the quotient of the completed free algebra in two variables $x, y$ by the ideal $[y,[x,y]]$.}

If $\omega_{k+a} \in \Filt^{k+a-\delta} \wedge^{k+a}\LL_A$, then 
$\xi^k\omega_{k+a} \in \Filt^{a-\delta}\wedge^a\LL_A$, which is zero 
if $ a < \delta$, giving (ii).
For (iii), observe the triangle 
$$ F^p\LL\Omega_A @>>>  \LL\Omega_A @>>>  \LL\Omega_A/F^p\LL\Omega_A @>{+1}>> $$
gives that $\omega $ is the transgression of a class $\nu_0 + \nu_1 + \dots 
+ \nu_{p-1} \in H^{i-1}(\Spec A, \LL\Omega_A/F^p\LL\Omega_A )$,
so $\omega = d \nu_{p-1}$, where $\nu_{p-1} \in (\wedge^{p-1}\LL_A)^{i-p}$ and
$D \nu_{p-1} = d \nu_{p-2} $, % \in (\wedge^{p-1}\LL_A)^{i-p+1}$, 
so $D\omega = 0$.

Hence $e^\xi \omega = \omega + \xi \omega + \dots + 
\frac{\xi^{p-\delta}}{(p-\delta)!}\omega$, a finite sum. 

\enddemo

\remark{Remark} If $\omega \in F^p\LL\Omega_A$ is transgressed
from a class in $\Filt^\delta (\LL\Omega_A/F^p\LL\Omega_A)$,
then $\omega \in \Filt^{\delta+1}F^p\LL\Omega_A$. In particular,
its class in $\gr^\delta F^p\LL\Omega_A$ is zero.
\endremark

\remark{Remark} The action of vector fields $\LL_A^*$ on $\wedge^i \LL_A$ is
well defined (independant of the choice of a cofibrant replacement for $A$).
The independence of $e^\xi$ on choices is a different phenomona, 
%. {\it is it? is $e^{D_B\nu}\omega $ equal to $\omega$ in the de Rham complex?}
as $e^{\xi + D_B\nu}\omega $ need not equal $e^\xi\omega$ in the de Rham complex.
%,if $d\nu \neq 0$.
%, as our notationis slightly imprecise. 
To define $e^\xi$,
we choose a cofibrant replacement for $B$, and a representative $\xi \in \Hom_B(M[-1],B)^0$ for the class of $\xi$. This defines a cofibrant 
representative for $A = \Sym^\xi_BM$, an $E \in (\LL_A)^0$ with $DE = \xi$,
and the morphism $e^\xi$.
Different cofibrant choices for representatives of $B$, $M$ and $\xi$ produce
canonically quis results.
%the independance of the above on choice of  $\xi$ is more subtle or not?
So our notation is harmlessly imprecise.
\endremark

\remark{Remark}
We only use part (iii) of the lemma, but the following may be clarifying.
The morphism $e^\xi : \LL\Omega_{\hat{A}}\to \LL\Omega_\bAh$ is a quis,
$\hat A =  \hat B \otimes_B A$. If $H^i(M) = 0 $ 
for $i \geq -1$, then $\LL\Omega_{B} \to \LL\Omega_{\hat A}$,
$\LL\Omega_{B} \to \LL\Omega_\bAh$ are quis; regardless 
%$H^i(\Spec \hat A, \LL\Omega_{\hat A})$ equals 
%$H^i(\Spec H^0(\hat A), \LL\Omega_{H^0(\hat A)})$,
$H^i(\Spec \bAh, \LL\Omega_{\bAh})$ equals 
$H^i(\Spec H^0(\hat A), \LL\Omega_{H^0(\hat A)})$,
and this coincides with Hartshorne-Deligne's algebraic deRham cohomology
of $H^0(A)$.
\endremark

\medskip

\head 3. A Darboux theorem % Endgame 
\endhead

Throughout this section we fix an invertible element $\CP$ of $\ch( k)$,
so $\CP = k[-d]$, % or $\CP = B[-d] \otimes \oddline$, 
a sign $\lambda_\CP \in \{\pm 1\}$,
and we assume 
$d \geq 1$. Given $B \in \cdgal$, and $ M\in B\mod$, 
write $\CP^{\pm 1}$ for $\CP^{\pm 1}\otimes_k B$, and
define 
$M^\dagger = \Hom_B(M,\CP^{-1}) \in B\mod$. Note that if
 $M$ is cofibrant and perfect then so is $M^\dagger$, 
and $M \quis M^{\dagger\dagger}$.

For $A \in \cdga$, $\Perf A$ becomes an exact category with 
weak equivalences and duality \cite{Sc,7.4},
where we define $\eta_M: M \to M^{\dagger\dagger}$ to be the natural
duality quis \cite{Sc,6.1} composed with $\lambda_\CP$.

Given integers $a \leq b $ with $a+b = -d$, consider the full subcategory of 
$\Perf A$ consisting of complexes with Tor amplitude in $[a,b]$;
this is the subcategory of complexes
%weakly equivalent
quis to one of the form  $M^a \to \cdots \to M^b$ 
with 
each $M_i[i]$ a summand of a free $A$-module $A^{n_i}$. This subcategory
inherits the structure of an exact category with weak equivalences and
duality.
%and if $a =b$ or $M$ is cofibrant, and we use the naive definition of 
% $\Hom_B(M,\CP)$,
%the duality weak equivalence $\eta_M : M \to M^{\dagger\dagger}$ is 
%an isomorphism as well as a quis.

%\footnote{The cofibrant objects are also semi-idempotent complete, as fibrations are degreewise surjections.}

\medskip

We remark that
%, having fixed either of the two possible natural choices of duality,  
there is a unique choice of $\lambda_\CP$ for which
$\CP$-shifted symplectic structures exist, and we may as well fix it.

\medskip
Moreover, we may also fix a cofibrant finitely presented $R \in \cdgal$,
and assume all cdga's lie over $R$,  $A \in \cdga_{R\backslash \cdot}$. Then
all of the theorems below are true, where we interpret $\LL_A$ to be
$\LL_{A/R}$ etc.

\sec{3.1} 
A {\it $\CP$-symmetric complex} on $A$ is a pair $M \in \Perf A$,
and $\varphi: M^\dagger \quis M$ a quis such that $\varphi^\dagger  = \eta_M\circ \varphi : M^\dagger \to M \quis M^{\dagger\dagger}$.

\smallskip

If $M$ is a  $\CP$-symmetric complex on $A$, and $N \hookrightarrow M$ is
a cofibration in $\Perf A$ equipped with a 
factorisation $(M/N)^\dagger \to N \to M$ of 
$(M/N)^\dagger \hookrightarrow M$,
we say that $N \to M$ is {\it co-isotropic}. For co-isotropic
$N \to M$,  
the quotient $ N/(M/N)^\dagger $ is a  $\CP$-symmetric complex on $A$.
When this quotient is quis to zero, we say that
$N \to M$,  or the triangle $N \to M \to N^\dagger$, is {\it Lagrangian}.

\smallskip

If $M$ is a $\CP$-symmetric complex on $A$, then the class of $M$
is {\it metabolic}
% {\it zero in the Witt group} of $\CP$-symmetric forms
if there is a Lagrangian 
$N \hookrightarrow M$; and it is
 {\it zero in the Witt group} of $\CP$-symmetric forms
if there is a metabolic complex $X$ such that $X \oplus M $ is metabolic.

\proclaim{Proposition} A $\CP$-symmetric complex $M$ is zero 
in the Witt group if and only if it is metabolic. Moreover, if this is so
there is a Lagrangian $S \to M \to S^\dagger$ with $H^i(S^\dagger) = 0$
for $i \geq - \intpart{ \frac{d-1}{2} }$.
\endproclaim
\demo{Proof} This is a consequence of `algebraic surgery'; see \cite{Sc,\S6}
for example, for a careful proof. Here is a sketch of 
the first statement: 
If the class of $M$ is zero in the Witt  
group, we have a triangle $L \to \cone(P^\dagger[-1] \to P) \oplus M \to L^\dagger @>{+1}>> $ for some $P \in \Perf A$ and morphism $P^\dagger[-1] \to P$,
hence we have a triangle $\tilde{L}\to M\to \tilde{L}^\dagger @>{+1}>> $
where $\tilde{L} = \cone(L \to P^\dagger)[-1]$.
\enddemo

\sec{3.2}
For any $M \in \Perf B$, the inclusion $M \hookrightarrow \Sym_B M = \bA$ 
induces $B \to \Sym_B M \otimes_B M^*$ by adjunction, and hence 
$\bA \to \bA \otimes_B M^*$. When $M = \LL_B^\dagger$, 
so $M^* = (\CP^{-1})^* \otimes \LL_B $,
the composite map $\bar\lambda : \bA \to \CP \otimes \LL_\bA$ 
is called the Liouville 
form, we have $D\bar\lambda = 0$, and $d\bar\lambda = (d+D)\bar\lambda 
= \bar\omega^\std$
is the ``standard'' shifted symplectic structure on $\bA$, which is 
manifestly deRham closed and non-degenerate.

Moreover, $\bar\omega^\std : \LL_\bA^\dagger \to \LL_\bA$ is a $\CP$-symmetric
complex, zero in the Witt group of $\bA$, as we have the vertical maps
in the diagram below
are quis's induced by $\bar\omega^\std$.
$$\CD
\LL_B\otimes_B\bA @>>> &  \LL_\bA @>>> & \LL_B^\dagger \otimes_B \bA @>{+1}>> \\
@AAA &  @AAA &  @AAA \\
\LL_B^{\dagger\dagger}\otimes_B\bA @>>> & \LL_\bA^\dagger @>>> & \LL_B^\dagger \otimes_B \bA @>{+1}>> \\
\endCD$$

Now let  $\txi \in H^{2}(\Spec B, \CP \otimes F^1\LL\Omega_B)$ be a one form,
and write $[\txi]_1 = \txi_1 \in H^1(\CP \otimes \LL_B)$. Regard 
$\txi_1$ as an element of $\CP[1]\otimes\LL_\bA$ via 
$\LL_B\otimes_B \bA \to \LL_\bA$. 
Then there exists an element $\xi' \in H^1(\CP \otimes \LL_\bA^\dagger)$
which is the image of a class $\xi \in \CP[1] \otimes \LL_B^{\dagger\dagger}
= \CP \otimes \Hom_B(\LL_B^\dagger,\CP^{-1}[1]) = \Hom_B(\LL_B^\dagger,B[1])$ 
such that $\xi'\bar\omega^\std = \txi_1$,
where $\xi' = \xi\otimes 1 \in \LL_\bA^\dagger \otimes \CP[1] = \LL_\bA^*[1]$.

Now $\xi : \LL_B^\dagger \to B[1]$ defines a deformation 
$A=\Sym^\xi_B\LL_B^\dagger$;
writing $\lambda$, $\omega^\std$ for the deformations of $\bar\lambda$, 
$\bar\omega^\std$ to $A$, we get $D_A \lambda = \Lie_\xi \lambda = \txi_1 
\in (\CP\otimes \LL_A)^1$, and $d\lambda = \omega^\std$,
so $(d+D_A)\lambda = \txi_1 + \omega^\std$.

Let $\omega^\std_\txi =  - \txi + (d+D_A)\lambda 
\in \LL\Omega_A$. This defines a class in 
$H^{2}(\Spec A, \CP\otimes F^2\LL\Omega_A)$
with $[\omega^\std_\txi]_2 = - \txi_2 + \bar\omega^\std$.
As $\gr [\omega^\std_\txi]_2 = \gr \,\bar\omega^\std$ is the identity map
in $\LL_B \otimes \LL_B^*\otimes_B A$, 
$[\omega^\std_\txi]_2 $ is non-degenerate, and $\omega^\std_\txi$ 
defines a twisted symmetric structure on $A = \Sym_B^\xi \LL_B^\dagger$, 
which we call the ``standard'' twisted structure attached to $\txi$.

\proclaim{\lsec{3.3} Proposition}
Let $A \in \cdgal$ have $H^0(A)$ f.p.{}, and let $\omega_2 : \LL_A^\dagger
\to \LL_A$ be a  $\CP$-symmetric complex.

If the class of $\omega_2$ in the Witt group of $\CP$-symmetric forms is 
zero, then there exists $B\to A$  in $\cdgal$ such that
$\LL_B^\dagger \otimes_B A \quis \LL_{A/B}$, $B$ is f.p.{}, 
and $H^i(\LL_{A/B}) = 0 $ for $i \geq - \intpart{ \frac{d-1}{2} }$.
\endproclaim
\demo{Proof} 
By proposition 3.1 we have a Lagrangian 
$S \to \LL_A \to S^\dagger @>{+1}>>$
with $H^i(S^\dagger) = 0$ for $i \geq -\intpart{\frac{d-1}2}$.
 The Tor amplitude of $S^\dagger$ is
in $[r, -\intpart{ \frac{d+1}{2}} ]$ for some $r$,
hence the Tor amplitude of $S$ is 
in $ [  - \intpart{ \frac{d+1}{2} }, 0]$, and so $\tordim S \leq  \intpart{ \frac{d}{2} }$, and $r = -d$.
The result is now immediate from the `Frobenius theorem' proposition 1.4.
\comment

As $\tordim \LL_A = d$, the Tor amplitude of $S^\dagger$ is contained in 
$[-d, -\intpart{\frac{d-1}2}]$, and that of $S$ is contained in 
$[-\intpart{\frac{d}2},0]$.

The result now follows from the `Frobenius theorem' proposition 1.?,
\footnote{It is easy to prove that the connectivity of $S^\dagger$ 
allows one to deduce the claim of proposition 1.? directly from proposition 1.1.
These extra few words about factoring $d$ hopefully make it clearer why such 
an argument works.}
if we can factor $\wedge^2\beta \circ d : \LL_A \to \wedge^2 S^\dagger$
through $\beta : \LL_A \to S^\dagger$. Such a map exists precisely when
$\wedge^2\beta\circ d\circ \alpha$ is homotopic to zero.

But $S \quis (S_{-r} \to \dots \to S_{0})$ where each $S_i[-i]$ is a 
summand of a free module $A^{n_i}$, so an element of $S$ is of the form
$\sum_{0\leq i \leq r} a_i v_i$, with $a_i \in A$, and $v_i \in S^{-i}$ of
degree $-i$. But $d\alpha(a_iv_i) = \pm a_i d \alpha( v_i)  + 
da_i.\alpha(v_i)$,  and $\wedge^2\beta ( d \alpha(v_i)) \in (\wedge^2 S^\dagger)^{-i}
= 0$, as $-i \geq -2 \intpart{\frac{d-1}2}$, so $\wedge^2\beta \circ d\circ \alpha = \beta d \wedge \beta \alpha : A \otimes S \to \wedge^2 S^\dagger$;
this is homotopic to zero as $\beta \alpha$ is by assumption.
\endcomment

\comment

\bigskip
{\it Proof by hand. But only part here}
By corollary 1.?, $\LL_A$ is perfect and $A$ is f.p.
By proposition 1.?, $A$ is a colimit of algebras $B_a$, where for $a \geq 1$
$B_{a+1}$ is built out of $B_a$ by attaching $-a$-cells.

By proposition 3.1 we have a Lagrangian $S @>>> \LL_A @>>> S^\dagger @>{+1}>>$
with $H^i(S^\dagger) = 0$ for $i \geq -\intpart{\frac{d-1}2}$.
Hence the map $\LL_{B_i}\otimes_{B_i}A \to S^\dagger$ is the zero map
if $i -1 \leq \intpart{\frac{d-1}2}$, and in particular 
for $i = \intpart{\frac{d}2}$. Set $r =  \intpart{\frac{d}2}$,
and choose a morphism $\LL_{B_r}\otimes_{B_r}A \to S$ factoring
$\LL_{B_r}\otimes_{B_r}A \to \LL_A$.

Consider the image of $H^{-r}(S)$ in  $H^{-r}(S/\LL_{B_r}\otimes_{B_r}A)$;
this is a finitely generated $H^0(A)$ module.
Let $w_1,\dots,w_\alpha$ be elements which generate it;
these define elements in $H^{-r}(\LL_{A/B_r})$,
which, as in the proof of prop 1.?, transgress to elements of $H^{-r+1}(B_r)$.
Choose lifts $\gamma'_1,\dots,\gamma'_\alpha \in  (B_r)^{-r+1}$. 

Set $B= B_r[z'_1,\dots,z'_\alpha \mid  D z'_i = \gamma'_i ]$, $\deg z'_i = -r$,
and define a map $B \to A$ by sending $z'_i$ to some choice of 
an element $z_i \in A^{-r}$ such that $Dz_i $ equals the image of $\gamma'_i$
in $A^{-r+1}$.

We can find a map 
$\LL_B \otimes_B A \to S$ %be a map 
factoring  the map
$\LL_B \otimes_B A \to \LL_A$,  extending
$\LL_{B_r} \otimes_{B_r} A \to \LL_A$ and sending $dz'_i$ to 
$w_i \in H^r(S/\LL_{B_r}\otimes_{B_r} A)$. % Writing $S/\LL_B$ for the cone of
So we have  $H^{-r}(S/\LL_B\otimes_B A) = 0$ by construction, and

\endcomment
%By construction 
%It is clear that $B \to A$ has the properties claimed.

\enddemo

\comment
\remark{Remark} It is quite easy to prove that this proposition holds locally
on $A$. So one could also prove it by showing lifts are unique-ish
 We were surprised that this stronger result is true.
\endremark
\endcomment

%\proclaim{Lemma} Let $B\to A$ be in $\cdgal$ such that $\LL_B^\dagger \otimes_B A \quis \LL_{A/B}$, $B$ is f.p.{}, and $H^i(\LL_{A/B}) = 0 $ for $i \geq - \intpart{ \frac{d-1}{2} }$. Then there exists $\xi \in H^0\Hom_B(\LL_B^\dagger, B[1])$ and a quis $\Sym_B^\xi \LL_B^\dagger \quis A$. \endproclaim

\proclaim{Lemma} Let $B\to A$ be in $\cdgal$ such that
$\bom_2 : \LL_B^\dagger \otimes_B A \quis \LL_{A/B}$, $B$ is f.p.{}, 
and $H^i(\LL_{A/B}) = 0 $ for $i \geq - \intpart{ \frac{d-1}{2} }$.
Then there exists $\xi \in H^0\Hom_B(\LL_B^\dagger, B[1])$ and a quis
$\alpha : \Sym_B^\xi \LL_B^\dagger \quis A$.

Moreover,
if $[d, \bom_2] = 0$, then we can choose this quis so that
$\alpha^*\bom_2 : \LL_B^\dagger \otimes_B \Sym_B^\xi\LL_B^\dagger 
\to \LL_B^\dagger \otimes_B \Sym_B^\xi\LL_B^\dagger$
is the pullback of a morphism $\LL_B^\dagger \to \LL_B^\dagger$.
\endproclaim
\demo{Proof} Arguing as in the previous proposition, 
the Tor amplitude of $\LL_{A/B}$ is 
%in $[r, -\intpart{ \frac{d+1}{2}} ]$ for some $r$,
%hence the Tor amplitude of $\LL_B$ and of $\LL_{A/B}^\dagger $ is 
%in $ [  - \intpart{ \frac{d+1}{2} }, 0]$, and so $\tordim B \leq  \intpart{ \frac{d+1}{2} }$, and $r = -d$.
in  $ [-d,  - \intpart{ \frac{d+1}{2} }]$. 

Hence by proposition 1.4, we may assume $B \to A$ is cofibrant, with
$A = B[y_1,\dots,y_a \mid D y_i = h_i ]$, 
for some $h_i \in B[y_1,\dots y_{i-1}]$
with $  \intpart{ \frac{d+1}{2} } \leq - \deg y_i \leq d$ for all $i$.
Hence $-\deg y_i - \deg y_j  \geq 2  \intpart{ \frac{d+1}{2} } > d-1$,
and only terms at most linear in $y_i$ appear in $Dy_i$. Write 
$Dy_i = h_i = \lambda_i + \sum_{j < i} \mu_{ij}y_j$ 
for some $\lambda_i, \mu_{ij}\in B$. 

Set $M = B\langle y_1,\dots y_a \mid Dy_i = h_i - \lambda_i\rangle \in B\mod$,
and $\xi : M \to B[1]$, $\xi(y_i ) = \lambda_i$, so that $\Sym_B^\xi M \to A$
is a quis. As $1\otimes_B d : A\otimes_B M \to \LL_{A/B}$ is a quis, we have 
%$\LL_B^\dagger$ and $M$ are two $B$-lattices in $\LL_{A/B}$, hence 
there exists a 
quis $\LL_B^\dagger \quis M$ %, {\it Elaborate!} 
by lemma 1.6, proving the result.

Now write $A = \Sym^\xi_B \LL_B^\dagger$, and $\bom_2$ for the pullback of 
$\bom_2$ to $A$. The filtration by degree on $A=\Sym^\xi_B\LL_B^\dagger$ 
induces one on $\Hom_A(\LL_A,\LL_A)$, the $i$'th piece of which is 
endomorphisms which raise filtration degree by $i$;
a morphism is in filtration degree 0 precisely when it is the pullback of
an endomorphism in $\Hom_B(\LL_B^\dagger,\LL_B^\dagger)$.

Write $\bom_2(dy_i) = \sum a_{ij}dy_j$, with $\deg a_{ij} = 
-\deg y_j + \deg y_i \leq 0$. As $\intpart{\frac{d+1}2} \leq -\deg y_i \leq d$,
$a_{ij}\in B$ unless $d$ is even, $\deg y_i = - d$, 
$\deg a_{ij} = -d/2 = \deg y_j$.

Hence if $d$ is odd, 
%$\bom_2 \in \Hom(\LL_B^\dagger,\LL_B^\dagger)\otimes_BA$,
$\bom_2$ is in filtration degree 0,
and if $d$ is even we can write 
$\bom_2(dy_i) = \sum \gamma_{i,jk}\, y_j dy_k + {\text{ terms of smaller 
filtration degree}}$,
where $\gamma_{i,jk} \in B^0$ is zero unless
$1/2\deg y_i = \deg y_j = \deg y_k = -d/2$.

By assumption $d \bom_2 (dy_i) = 0 $ in $\wedge^2 \LL_{A/B}$, so 
$\sum_{j,k} \gamma_{i,jk}\, dy_j dy_k  = 0$, and hence 
 $\bom_2(dy_i) =d ( \sum_{j <k} \gamma_{i,jk}\, y_j y_k )$.
%{\it Check! Really its zero in cohom. Be more precise,}

Now define $\LL_B^\dagger \to A  $ by
$y_i \mapsto y_i - \sum_{j <k} \gamma_{i,jk}\, y_j y_k $,
and let $\alpha : A \to A$ be the induced morphism of cdgas.
Note that $D_B(y_j) = 0 $ if $\deg y_j = -d/2$, so $D\alpha = \alpha D$,
$\alpha $ is invertible, and $\alpha^*\bom_2 :
 \LL_B^\dagger \otimes_B A  \to  \LL_B^\dagger \otimes_B A $ is in filtration
degree zero.

\enddemo
\remark{Remark} Note that $\alpha^* \xi \neq \xi$, so the
isomorphism chosen in the second part of the lemma changes the choice 
of $\xi$ from the first part.
\endremark

\proclaim{\lsec{3.4} Lemma} Let $B \in \cdgal$ with $H^0(B)$ f.p.{}, $\xi \in H^1(\LL_B^\dagger)$, 
$A = \Sym_B^\xi \LL_B^\dagger$, and suppose given 
$\omega \in H^{2}(\Spec A, \CP\otimes F^2\LL\Omega_A)$ such that

i) $[\omega]_2 \in H^{0}(\CP\otimes\wedge^2\LL_A)$ defines a quis
$\LL_B \otimes_B A \quis \LL_{A/B} \quis \LL_B^\dagger\otimes_B A$, and 

ii) the class of $\omega$ in $H^{2}(\Spec A, \CP\otimes\LL\Omega_A)$ is zero.

\noindent
Then there exists an $f : B \to \CP[1]$ such that 
$$e^\xi\omega = \pi^*(df) \in H^2(\Spec\bA, \CP \otimes F^1\LL\Omega_\bA),$$
where $\pi: B \hookrightarrow \bA = \Sym_B\LL_B^\dagger$ is the inclusion.

%Then there exists a $\tilde{\xi} \in H^{1}(\Spec B,\CP\otimes F^1 \LL\Omega_B)$
%such that $[\tilde{\xi}]_1 = \tilde{\xi}_1$, 
%where $\tilde{\xi}_1 \in H^1(\LL_B)$ is the one form defined by $\xi$.

%Moreover, if $H^{1}(\Spec B,\CP\otimes \LL\Omega_B) = 0$, then $\tilde{\xi}$ is unique.
\endproclaim

\remark{Remark} We have that the image of $H^i(\Spec A,\CP\otimes F^2\LL\Omega_A) \to H^i(\Spec A,\CP\otimes \LL\Omega_A)$ is zero for $ i \leq  d$,
so the condition (ii) in the lemma is always satisfied if $d > 1$.
\endremark

\demo{Proof} Choose a cofibrant replacement for $B$, and 
$\xi \in (\LL_B^\dagger)^1$ with class $\xi$, and consider 
$e^\xi : \LL\Omega_A \to \LL\Omega_\bAh$ as in lemma 2.7.

We have $\omega_2 \in \Filt^1\wedge^2\LL_A$ by hypothesis, and
by hypothesis $\omega$ is transgressed from a class in $\LL\Omega_A/F^2\LL\Omega_A$.
So, by
lemma 2.7(iii) we may assume that $\omega = \omega_2$
and $e^\xi\omega = \xi \omega_2 + \omega_2 
\in H^2(\Spec \bA, \CP \otimes \Filt^1 F^1 \LL\Omega_\bA)$.

If we write $ (\xi \omega_2 + \omega_2 )^\grade{\lambda}$
for the $\lambda$'th graded piece of $\xi \omega_2 + \omega_2 $ as in \S2.3,
then for $\lambda > 0$, $(\xi \omega_2 + \omega_2)^\grade{\lambda} 
= (d+\bD) D_\Kos( (\xi \omega_2 + \omega_2 )^\grade{\lambda}/\lambda)
= (d+\bD) D_\Kos (\omega_2^\grade{\lambda})$,
as $D_\Kos (\xi\omega_2) = 0$, and this is zero in $H^2(\Spec \bA, \CP\otimes F^1\LL\Omega_\bA)$.

Hence $\xi \omega_2 + \omega_2  = (\xi \omega_2 + \omega_2 )^\grade{0} 
= \pi^* \gamma$, for $\gamma \in H^2(\Spec B, \CP\otimes F^1\LL\Omega_B)$.
Moreover, the class of $\gamma$ in $H^2(\Spec B, \CP\otimes \LL\Omega_B)$ is zero. This is because the class of $e^\xi\omega $ is zero in 
$H^2(\Spec \bA, \CP\otimes \LL\Omega_\bA)$ by hypothesis (ii),
and $\pi^*: H^i(\Spec B,\CP \otimes \LL\Omega_B) \to H^i(\Spec \bA,\CP \otimes \LL\Omega_\bA)$ is an isomorphism for all $i$.
%This is clear if $d \geq 2$, as $H^i(\Spec B,\CP \otimes \LL\Omega_B) = 0$ for $i < d$, and for $d =1 $ it follows as there is no constant term in $\gamma$.

Hence $\gamma = df$, for some element $f \in H^1(\Spec B,\CP)$.

%regard $\omega \in H^{2-d}(\Spec A, \Filt^1F^2\LL\Omega_A)$.

\enddemo

\sec{3.5}
Let $B$, $\xi$, $A$, $\omega$ be as in the lemma and its proof, 
and suppose also that
$[\omega]_2 : \LL_B^\dagger \otimes_B A \to \LL_B^\dagger \otimes_B A $
is the pullback of a map $\sigma : \LL_B^\dagger \to \LL_B^\dagger $.
Define $A' = \Sym_B^{\xi\sigma} \LL_B^\dagger$, and 
$\sigma : A' =  \Sym_B^{\xi\sigma} \LL_B^\dagger \to A = \Sym_B^\xi \LL_B^\dagger$
to be the induced map in $\cdga_{B\backslash\cdot}$. Observe that the 
%induced 
map 
$\LL_B^\dagger \otimes_B A \quis \LL_{A'/B}\otimes_{A'}A @>{d\sigma}>>
\LL_{A/B} \quis \LL_B^\dagger \otimes_B A $
is %just 
induced by $[\omega]_2$ and so is a quis, i.e.{} $\LL_{A'/A} = 0$.
Moreover, as $H^0(B) \twoheadrightarrow H^0(A')$, $H^0(B) \twoheadrightarrow H^0(A)$ are surjections, so is $H^0(A')\twoheadrightarrow H^0(A)$, hence $A'\to A$ is a quis,
by corollary 1.1(ii).

Put $\teta = (\xi\omega_2 + \omega_2)^\grade{0} = \xi\omega_2 + \omega_2^\grade{0}$. Then $\eta = \xi \sigma$, and $\sigma^* \omega^\std_\teta = 
\sigma(\omega^\std_\teta) \in \wedge^2\LL_A$ satisfies
$\sigma^*\omega^\std_\teta - \omega \in \Filt^0\wedge^2\LL_A$.

%As $\teta = \pi^*(df)$ for an element $f \in H^1(\Spec B,\CP)$,
%we have $ \omega_2^\grade{0} = 0 $ in $H^1(\Spec A', \CP \otimes F^1\LL\Omega_{A'})$.

\proclaim{Lemma} $\sigma^* \omega^\std_\teta = \omega $ in
$H^2(\Spec A, \CP \otimes F^2\LL\Omega_A)$
\endproclaim
\demo{Proof} Put $\delta = \sigma^* \omega^\std_\teta - \omega $. Consider
the filtration by degree, $\bFilt$. If $\delta $ is in $\bFilt^\lambda
\wedge^2\LL_A$ 
and $\lambda > 0$, then
$\bgr^\lambda(\delta) = (d+D) D_\Kos (\bgr^\lambda(\delta)/\lambda)$.
But $D_\Kos(\bgr^\lambda\delta) = 0$, as $\delta \in \Filt^0\wedge^2\LL_A$.
So $\delta$ is in $\bFilt^0$ and so pulled back from $B$.
%Finally, observe 
It follows 
that by construction $\delta$ is 
%a multiple of $\sigma^* \omega_2^\grade{0}$, which is 
zero.
\enddemo

\sec{3.6} To summarise, it seems we have proved the following.

\proclaim{Theorem} Let $\CX \to \Cal R$ be a connected Deligne-Mumford dg-stack
over some base DM dg-stack $\Cal R$, with 
% with perfect cotangent complex, 
$\CP$ an invertible complex of %$\bz/2\bz$-graded 
$D$-modules on $\CX$,
$H^{d}(\CP) \neq 0$ for some $ d > 0$, 
and $\omega \in H^0(\CX,H^2(\CP\otimes F^2\LL\Omega_{\CX/\Cal R}))$
a $\CP$-shifted symplectic form.

\noindent
Suppose that $\CD\Spec A @>{i}>> \CX \\ @VVV @VVV \\ \Spec R @>>> \Cal R \endCD$
is an etale map for which 

i) the class of $i^*[\omega]_2$ is zero in the Witt group of $A$,

ii) the class of $i^*\omega$ is  zero in 
$H^2(\Spec A, \CP \otimes \LL\Omega_{A/R})$, and

iii) the underlying $D$-module of $i^*\CP$ is trivial.\footnote{This assumption is somewhat harsh. We will repost this note with a more general statement later.}

\noindent
Then there is a f.p.{} cdga $B$,  $f : B \to \CP[1]$,  and quis 
$\sigma : A' = \Sym_B^{\eta} \LL_B^\dagger \to A$ such that
$\sigma^*\omega^\std_{df} = i^* \omega$, where $\eta : \LL_B^\dagger \to B[1]$
is induced from $df$.

\endproclaim

\demo{Proof} Given $A$ as above, 
%\footnote{and writing $k$ for $R$ as we do in the rest of the paper.}
 observe that by lemma 1.7 the underlying $A$-module of $i^*\CP$ is $\CO[-d]$, necessarily with $d \geq 0$. As the class of 
$[\omega]_2$ in the Witt group is zero, by proposition 3.3 and lemma 3.3(i)
we have a finitely presented cdga  $ B \to A$, and $\xi : \LL_B^\dagger \to B[1]$
such that $\Sym^\xi_B \LL_B^\dagger \to A$ is a quis. 
As the class of $i^*\omega$ is zero in the deRham complex of $A$, we can 
assume that $i^*\omega$ is transgressed from $\LL_A/F^2\LL_A$,
so by lemma 3.3(ii) we can further assume $\bom_2$ is pulled back from
a morphism $\LL_B^\dagger \to \LL_B^\dagger$. Lemma 3.4 now gives an
$f : B \to \CP[1]$ with $e^\xi\omega = \pi^*df$, and the discussion 
and lemma in \S3.5 gives the result.

\enddemo

\remark{Remark} If $d \neq 2 \Mod 4$, every geometric 
point $x \in X$ has a %n etale 
neighbourhood for which (i) holds.We can do slightly better.
Let $\CX$ be a Deligne-Mumford dg-stack, $\omega$ a $\CP$-shifted symplectic
form on $\CX$, $\Spec A' \to \CX$ an etale map, and $x \in \CX$ a closed point.

%As minimal resolutions over (dg) local rings are unique, a
Algebraic surgery
and Witt cancellation 
gives the 
localisation $(\LL_{A'})_{(x)}$ is quis to a direct sum of perfect complexes
$P_{(x)} \oplus M_{(x)}$, where $M_{(x)}$ is metabolic and $P_{(x)}$ has
Tor amplitude in $[-d/2,-d/2]$; moreover we can insist that $P_{(x)}$ is
zero if $d \neq 2 \Mod 4$. As $\CX$ is locally finitely presented,
and $M_{(x)}$, $P_{(x)}$ are perfect, there is some Zariski neighbourhood
$A'\to A$ of $x$, and $P, M \in \Perf A$, such that 
$\LL_A = P \oplus M$, the localisations of $P, M$ are $P_{(x)}, M_{(x)}$,
and $M$ is metabolic and $P$ has Tor amplitude $[-d/2,-d/2]$.\footnote{
By making an etale cover we may further insist that, if non-zero,
$P$ is just $A[d/2]$ 
with the quadratic form $x^2$.} %if we so choose.

Put $C = \Sym_A^\gamma P$, where $\gamma : P \to P$ is the identity map.
Then $M\otimes_A C \to \LL_C$ is a quis, so $\Spec C$ carries a $\CP$-shifted
symmetric structure with $\LL_C$ zero in the Witt group. Hence assuming
(ii) and (iii) of the theorem, we get that $C$ is a twisted shifted cotangent 
bundle. But $A \to \Sym_C(P)$ is a quis. We have proved

\endremark

\proclaim{Corollary} Let $x \in \CX$, $\CP = \CO[-d]$.
Then there is a neighbourhood $\CD \Spec A @>{i}>> \CX \\ @VVV @VVV \\
\Spec R @>>> \Cal R \endCD$ of $x$, and quis
 $\sigma:  \Sym^\eta_B(\LL_B^\dagger + P) \to A$, where 
$B \in \cdga_{R\backslash \cdot}$ is f.p.{}, $P \in \Perf B$ is a $\CP$-shifted
symmetric complex of Tor amplitude $[-d/2,-d/2]$, zero unless $d = 2 \Mod 4$,
and $f : B \to \CP[1]$, such that $i^*\omega$ is the pullback by $\sigma$ of
the  sum  of the 
standard symplectic form $\omega^\std_{df}$ with the form induced from $P$.
\endproclaim

\sec{3.7} It seems that the theorem extends to Artin dg-stacks with little extra effort: if $\CX \to \Cal R$ is an Artin dg-stack, $\omega$ a $\CP$-shifted symplectic structure, and $d > 0$, then locally $\CX$ is a twisted shifted cotangent bundle $\Sym^\eta_{\CO_{\CY}}(\LL_{\CY}^\dagger + P)$, where $\CY$ is an Artin dg-stack,  $H^i(\LL_{\CY}^\dagger) = 0$ for $ i \geq 0$, and $P \in \Perf \,\CY$ 
 is a $\CP$-shifted
symmetric complex of Tor amplitude $[-d/2,-d/2]$, zero unless $d = 2 \Mod 4$, as above. We will repost with details shortly.

\comment
$\alpha : L \to \CX$ is isotropic if $\CX$ is $\CP$-symplectic,
and $\nu \in H^2(L,\CP[1]\otimes F^2 \LL\Omega_\CX)$ has 
$\alpha^* \omega = \ (d+D) \nu$.

\endcomment

\comment

\smallskip
Observe that $[\bE, \bar\lambda] \in \wedge^2 \LL_\bA^*$ is the 
Poisson structure dual to $\bar\omega^\std$.

\bigskip

\sec{3.7}
We record belatedly\footnote{We do not use  these formulae, but the reader may find them reassuring, especially if they do not mind sign errors.}
that if we choose a cofibrant replacement for $B$, 
$B = k[z_1,\dots,z_n \mid Dz_i = f_i]$, $ f_i \in k[z_1,\dots,z_{i-1}]$,
then $\LL_B = B\langle dz_1,\dots,dz_n \mid D(dz_i)= - \sum_{j < i} 
\frac{\partial f_i}{\partial z_j} dz_j \rangle$, and if we set 
$y_i = (dz_i)^\dagger$, then 
$\bA = B[y_1,\dots,y_n \mid D_\bA(y_i) = - \sum_{j >  i}
\frac{\partial f_j}{\partial z_i} y_j ]$,
$\LL_\bA = A\langle dz_1, \dots, dz_n,  dy_1,\dots,dy_n \mid 
D_\bA(dz_i) =  - \sum_{j < i} \frac{\partial f_i}{\partial z_j} dz_j,
D_\bA(dy_i) = \sum_{j >  i} \frac{\partial f_j}{\partial z_i} dy_j +
\sum \frac{\partial^2f}{\partial z_k \partial z_i} dz_k . y_j
\rangle
$.
We have $\bar\lambda = \sum y_i dz_i$,
$\xi = \sum \xi_i (dy_i)^*$, $\tilde{\xi}_1 = \sum \xi_i dz_i$, for some $\xi_i \in B$.

%\smallskip
%$\lambda, E, etc$ on $A$

\endcomment

%We suppose given $A \in \cdgal$, and 

%\head 3. Appendix to \S2 \endhead

\enddocument